\documentclass[final,5p,times]{elsarticle}
\biboptions{numbers,sort&compress}
\usepackage{amsmath,amssymb,bbm,array,amsfonts,float,mathtools,amsthm}
\usepackage{mathrsfs}
\usepackage{booktabs}%
\usepackage{algorithm}%
\usepackage{algorithmicx}%
\usepackage[noend]{algpseudocode}%
\usepackage{listings}%
\usepackage{xcolor}
\usepackage{todonotes}
\usepackage{comment}
\usepackage[hyphens]{url}
\usepackage[hidelinks]{hyperref}
\hypersetup{breaklinks=true}

\newcommand{\MAE}{{\text MAE}}

\newcommand{\rG}{{\rm G}}

\newcommand{\rI}{{\rm I}}

\newcommand{\bw}{{\bf w}}

\newcommand{\cJ}{\mathcal{J}}

\newcommand{\sX}{\mathscr{X}}

\newcommand{\R}{\mathbb{R}}
\newcommand{\C}{\mathbb{C}}

\newcommand{\Spin}{{\rm Spin}}
\newcommand{\SU}{{\rm SU}}

\newcommand{\U}{{\rm U}}

\newcommand{\End}{{\mathrm{End}}}

\renewcommand{\epsilon}{\varepsilon}

\newcommand{\Hol}{\mathrm{Hol}}

\renewcommand{\Im}{\mathop{\mathrm{Im}}}

\renewcommand{\Re}{\mathop{\mathrm{Re}}}

\newcommand{\qandq}{\quad\text{and}\quad}

\newcommand{\qwithq}{\quad\text{with}\quad}

\def\<{\mathopen{}\left<}
\def\>{\right>\mathclose{}}
\def\({\mathopen{}\left(}
\def\){\right)\mathclose{}}

\newcommand{\floor}[1]{\lfloor#1\rfloor}

\usepackage{multicol, color}

\definecolor{gold}{rgb}{0.85,.66,0}
\definecolor{cherry}{rgb}{0.9,.1,.2}
\definecolor{burgundy}{rgb}{0.8,.2,.2}
\definecolor{orangered}{rgb}{0.85,.3,0}
\definecolor{orange}{rgb}{0.85,.4,0}
\definecolor{olive}{rgb}{.45,.4,0}
\definecolor{lime}{rgb}{.6,.9,0}
\definecolor{green}{rgb}{.2,.7,0}
\definecolor{grey}{rgb}{.4,.4,.2}
\definecolor{brown}{rgb}{.4,.3,.1}

\newtheorem{theorem}{Theorem}
\newtheorem{proposition}[theorem]{Proposition}

\newtheorem{conjecture}[theorem]{Conjecture}
\theoremstyle{definition}
\newtheorem{definition}[theorem]{Definition}

\journal{Physics Letters B}

\begin{document}

\begin{frontmatter}

\title{Machine Learning Sasakian and $\rG_2$ topology on contact Calabi-Yau $7$-manifolds}

\author[a]{Daattavya Aggarwal}
\author[b,c,d,e]{Yang-Hui He}
\author[b,c]{Elli Heyes}
\author[f]{Edward Hirst}
\author[g]{Henrique N. S\'a Earp}
\author[g]{Tom\'as S. R. Silva}

\affiliation[a]{
    organization={Department of Computer Science and Technology, University of Cambridge}, 
    postcode={CB3 0FD},
    country={UK}  
}
\affiliation[b]{
    organization={Department of Mathematics},
    addressline={City, University of London},
    postcode={EC1V 0HB},
    country={UK}
}
\affiliation[c]{
    organization={London Institute for Mathematical Sciences}, 
    addressline={Royal Institution}, 
    city={London},
    postcode={W1S 4BS},
    country={UK}
}
\affiliation[d]{
    organization={Merton College, University of Oxford}, 
    postcode={OX1 4JD},
    country={UK}
}
\affiliation[e]{
    organization={School of Physics, NanKai University}, 
    city={Tianjin}, 
    postcode={300071},
    country={P.\ R.\ China},
}
\affiliation[f]{
    organization={Centre for Theoretical Physics},
    addressline={Queen Mary, University of London},
    postcode={E1 4NS},
    country={UK}
}
\affiliation[g]{
    organization={Institute of Mathematics, Statistics and Scientific Computing, University of Campinas (Unicamp)}, 
    postcode={13083-859},
    country={Brazil}  
}

\begin{abstract}
We propose a machine learning approach to study topological quantities related to the Sasakian and $G_2$-geometries of contact Calabi-Yau $7$-manifolds. Specifically, we compute datasets for certain Sasakian Hodge numbers and for the Crowley-N\"ordstrom invariant of the natural $G_2$-structure of the $7$-dimensional link of a weighted projective Calabi-Yau $3$-fold hypersurface singularity, for 7549 of the 7555 possible $\mathbb{P}^4(\textbf{w})$ projective spaces. These topological quantities are then machine learnt with high performance scores, where learning the Sasakian Hodge numbers from the $\mathbb{P}^4(\textbf{w})$ weights alone, using both neural networks and a symbolic regressor which achieve $R^2$ scores of 0.969 and 0.993 respectively. Additionally, properties of the respective Gr\"obner bases are well-learnt, leading to a vast improvement in computation speeds which may be of independent interest. The data generation and analysis further induced novel conjectures to be raised. 
\end{abstract}

\begin{keyword}
$\rG_2$-manifolds, machine learning, Hodge numbers, Crowley-N\"ordstrom invariant, contact Calabi-Yau manifolds\\
\textit{Report Number}: QMUL-PH-23-14
\end{keyword}

\end{frontmatter}

\section{Introduction}
\label{sec:intro}


\paragraph{Motivation}

Since its inception in 2017 \cite{he2017,ruehle2017,carifio2017,krefl2017}, the study of Calabi-Yau (CY) manifolds in the context of string theory compactifications, using machine learning techniques has flourished, encompassing a wide array of investigations. Notably, these methods have been employed to predict Hodge numbers \cite{Bull:2019cij,He:2020lbz,Erbin_2021,Berman:2021mcw,pmlr-v197-aslan23a}, learn Ricci-flat Calabi-Yau metrics \cite{ashmore2020,anderson2020,jejjala2020,douglas2020,larfors2021}, forecast line bundle cohomologies \cite{klaewer2019}, generate new Calabi-Yau manifolds \cite{berglund2023new}, and uncover volume bounds on Sasaki-Einstein manifolds \cite{Manko:2022zfz}. Furthermore, machine learning techniques have found various other applications in geometry and physics \cite{He:2020eva,Bao:2020nbi,Bao:2021vxt,Bao:2021auj,Bao:2021olg,Bao:2021ofk,AriasTamargo:2022qgb,Dechant:2022ccf,Chen:2022jwd,Cheung:2022itk,Ashmore:2023ajy}. For an extensive review, see \cite{bao2022,he2023machine}.
As important as Calabi-Yau compactification is to string theory, $7$-manifolds of holonomy  $\rG_2$ are crucial to M-theory compactification  \cite{Acharya2001,Acharya2004}. Nonetheless, there have yet been no successful applications of machine learning in the context of $\rG_{2}$-geometry, let alone on compact manifolds with such special holonomy. This is no doubt due to the scarcity of dedicated databases for $\rG_2$-manifolds, which in turn reflects the difficulty in describing (torsion-free) $\rG_2$-manifolds systematically in terms of algebraic discrete data. 

While the holonomy reduction to $\rG_2$ amounts to a difficult non-linear PDE, pragmatically it may be relaxed in a number of ways, by considering the fundamental notion of a \emph{$\rG_2$-structure}: a non-degenerate  $3$-form $\varphi$ which induces a so-called \emph{$\rG_2$-metric}  $g_\varphi$;
its failure to give rise to a metric with $\Hol(g_\varphi)\subset\rG_2$  is encoded by its \emph{full torsion tensor} $T:=\nabla^{g_\varphi}\varphi$.
A $\rG_2$-structure is  \emph{closed} if  $d\varphi=0$ and \emph{coclosed} if $d\psi=0$, where $\psi:=\ast_\varphi\varphi$, and the torsion-free condition $T=0$ is equivalent to $\varphi$ being both closed and coclosed.  
We propose therefore to work with coclosed $\rG_2$-structures on certain \emph{contact Calabi-Yau (cCY)} $7$-manifolds, which are closely related to the weighted projective Calabi-Yau $3$-folds famously studied in \cite{CANDELAS1990383}.

Despite their unsuitability to M-theory, torsionful $\mathrm{G}_2$ structures retain  relevance in the context of (3+7)-dimensional  heterotic supergravity with flux, as demonstrated by \cite{delaOssa:2014lma,delaOssa:2017pqy,delaOssa:2017gjq}. Indeed, as shown by \cite{Lotay2023}, one can explicitly solve the corresponding Strominger system on cCY $7$-manifolds, by way of coclosed $\rm{G}_2$-structures together yielding non-trivial scalar and $\rm{G}_2$-instanton gauge fields, with constant dilaton, as well as an $H$-flux with prescribed Chern-Simons defect, in accordance to the `anomaly-free' condition referred to as the heterotic Bianchi identity.

\paragraph{Topological invariants of CY links}
Contact Calabi-Yau manifolds were introduced by Tomassini and Vezzoni in \cite{tomassini2007contact}, and they consist of Sasakian manifolds endowed with a closed basic complex volume form, which is `transversally holomorphic' in the sense of foliations. It was shown in \cite{Habib2015} that such a manifold carries naturally a coclosed $\rG_{2}$-structure. 

A special class of such structures arises from Calabi-Yau links, which were first discussed from the perspective of $\rG_2$ topology in \cite{Calvo-Andrade2020}. 
A $7$-dimensional weighted link $K_f$  is obtained as the intersection of a possibly small $S^{9}\subset \C^5$ with a weighted homogeneous affine variety (defined by the zero locus of the polynomial $f$) having an isolated singularity at the origin. Milnor showed that such links are $2$-connected compact smooth manifolds, indeed $K_f$ is the total space of a Hopf $S^1$-bundle over a (weighted) projective $3$-orbifold in $\mathbb{P}^4(\textbf{w})$, for appropriate choices of polynomial degree and weighted $\C^\times$-action, see \S\ref{sec:background}.
Interestingly, the dataset of possible weights that admit these CY $3$-folds consists of the 7555 cases classified in \cite{CANDELAS1990383}.
Therefore, we pursue the construction of a Calabi-Yau link for each of these weight systems, computing the following two types of  topological invariants.

From the perspective of Sasakian topology, the (basic) Hodge numbers $h^{p,q}$  can be obtained as the dimensions of certain linear subspaces of the Milnor algebra $\mathbb{M}_f=\C [[z_1,\dots,z_5]]/\cJ_f$, defined by the corresponding Jacobian ideal of $f$ \cite{Itoh2004}. We provide the first systematic computation of the Sasaki-Hodge numbers $\{h^{3,0},h^{2,1}\}$ for this class of $7$-dimensional CY links. 

On the other hand, considering their $\rG_2$-topology, a CY link bounds an $8$-dimensional Milnor fibre which smoothly extends the $\rG_2$-structure $\varphi$ as a spinor field, hence it is possible to explicitly compute the  Crowley-Nordstr\"om (CN)  homotopy invariant  $\nu(\varphi)\in \mathbb{Z}/{48}\mathbb{Z}$, introduced in \cite{CNInvariant}. Building upon the calculations first carried out in \cite{Calvo-Andrade2020}, we obtain an exhaustive dataset of $\nu$-invariants for Calabi-Yau links. 

\paragraph{Machine Learning cCY topology}

We analyse these two sets of topological data from a perspective similar to what has been done for Calabi-Yau manifolds \cite{he2017}.
In the standard Calabi-Yau case, the weights defining the ambient projective space are sufficient to uniquely determine the Calabi-Yau $3$-fold's Hodge numbers, motivating ML of the known formulas from weights to Hodge numbers \cite{Vafa:1989xc}. 
However, in the $7$-dimensional Calabi-Yau link case, no such explicit formula is known, and one would initially expect the specific polynomial coefficients chosen to change the topology.
Extending the ML techniques to these link invariants would establish existence of an approximate formula for Sasakian Hodge numbers from the weight information, from which ML interpretability techniques may be used to uncover its true form.
This formula would provide new insights into Sasakian structures, as well as being dramatically quicker to compute, as well as open the door for their application on other related invariants, including for $\rG_2$-structures as motivated by their ML in this work.

We therefore extend previous work learning CY Hodge numbers from weights, to predicting Calabi-Yau link topological properties (namely Sasakian Hodge numbers and CN invariants).
We find that, whilst the machine is able to learn the Sasakian Hodge number topology of these manifolds with high-performance measures, the same cannot be established for the CN invariant. 
The datasets of the weighted Calabi-Yau polynomials used in the link construction, with the computed Sasakian Hodge numbers and CN invariants, as well as the scripts used for analysis and machine learning, are available on \href{https://github.com/TomasSilva/MLcCY7.git}{\texttt{GitHub}} \cite{github}.

This letter is organised as follows: in \S\ref{sec:background} we survey some background to contact Calabi-Yau manifolds, $\rG_2$-geometry, and machine learning; in \S\ref{sec:method} we describe the methodology for the construction of the Calabi-Yau link data and perform relevant statistical analysis of the datasets of invariants; 
in \S\ref{sec:results} we present the results of the machine learning investigations; and we conclude in \S\ref{sec:conclusion}, discussing some future prospects.
\section{Background}
\label{sec:background}

\subsection{Calabi-Yau links}
\label{sec: CY links}

One may interpret structure group reductions on an odd-dimensional contact metric manifold $(K^{2n+1},\eta,\xi,g)$ as `even-dimensional' structures `transverse' with respect to a $S^1$-action along the fibres of a submersion $S^1\to K\to V$. Here $\eta\in\Omega^1(K)$ denotes the contact form and $\xi\in\sX(K)$ its  (unit) dual Reeb field, such that $\eta(\xi)=1$. Whenever clear from context, we will omit mention of the Riemannian metric $g$, for simplicity. 

In particular, Sasakian geometry may be seen as transverse Kähler geometry, corresponding to the reduction of the transverse holonomy group to $\U(n)$. These are equipped in addition with a transverse complex structure $J\in\End(TK)$ such that $J\circ J=-\rI_{TK}+\eta\otimes \xi$, yielding a decomposition of forms by basic bi-degree,  and a transverse symplectic form  $\omega=d\eta\in\Omega^{1,1}(K)$, all of which satisfy suitable compatibility conditions; for more details see eg. \cite[\S2]{Portilla2023} or the canonical reference \cite{Boyer2008}. Furthermore, Sasakian manifolds with special transverse holonomy $\SU(n)$ are studied by Habib and Vezzoni \cite[\S~6.2.1]{Habib2015}:
\begin{definition}
\label{def:cCY}
    A Sasakian manifold $(K^{2n+1},\eta,\xi,J,\Upsilon)$ is said to be a \emph{contact Calabi--Yau manifold} (cCY) if $\Upsilon$ is a nowhere-vanishing transverse form of horizontal type $(n,0)$, such that 
$$
\Upsilon\wedge\bar{\Upsilon}
=(-1)^{\frac{n(n+2)}{2}}\omega^n
\qandq 
d\Upsilon=0,
\quad\text{with}\quad
\omega=d\eta.
$$
\end{definition} 

A polynomial $f\in\C[z_1,\dots,z_{n+2}]$, for $n\geq 2$, is said to be \emph{weighted homogeneous} of degree $d$ with weight vector $\mathbf{w}=(w_1, \dots, w_{n+2})\in \mathbb{Z}_{>0}^{n+2}$,
if it is homogeneous of order $d$ with respect to the $\C^{\times}(\bw)$-action on $\C^{n+2}$ 
\[(t,z) \mapsto t\cdot z = (t^{w_1}z_1, \dots, t^{w_{n+2}}z_{n+2}).\]
Such an $f$ defines an affine variety 
\[V_f = (f)=\{z\in \C^{n+2} \mid f(z)=0\},
\] which, in general, admits a singularity at the origin.

Assuming that the origin is an isolated singularity, the intersection of $V_f$ with a surrounding small hypersphere $S_{\varepsilon}^{2n+3}$ is a compact smooth $(2n+1)$-manifold $K_f=V_f \cap S_{\varepsilon}^{2n+3}$, the so-called weighted link of the singularity \cite{Milnor1969}.
A weighted link $K_f$ of degree $d$ and weight $w$ is a \emph{Calabi-Yau link} if
\begin{equation}
\label{CYcondition}
    d = \sum_{i=1}^{n+2}w_i,   
\end{equation}
which precisely guarantees the existence of a cCY structure on $K_f$.
The dimension of the moduli space of these cCY structures is well-understood and discussed in \S\ref{sec:conclusion}.

\subsection{Sasakian Hodge numbers of a CY link}
\label{sec:sasakian-hodge}

The $\C^{\times}(\bw)$-action on $\C^{n+2}$  induces a contact-metric $S^1$-action on $K_f$. It admits finitely many distinct isotropy subgroups, contained in some finite subgroup $\Gamma \subset S^1$, so that $K_f$ admits a double fibration over a projective $n$-orbifold $V\subset\mathbb{P}^{n+1}(\mathbf{w})$,
\[
\pi: K_f \longrightarrow K_f/\Gamma \longrightarrow K_f/S^1 = (V_f \setminus \{0\})/ \C^{\times} \vcentcolon = V^*_f
\]
The following key theorem allows us to compute certain mixed Hodge numbers $h^{p,q}(K_f)$ from the dimensions of the primitive cohomology groups $H_0^{n}(V_f^{*})$, for $p+q=n$, which in turn can be obtained from the Milnor algebra. A brief survey of Sasakian Hodge numbers can be found in \ref{sec: Appendix}.

\begin{theorem}[{\cite[Theorem 1.2]{Itoh2004}, \cite{Steenbrink1977,Steenbrink1983}}]
\label{th:itoh}\mbox{}\\
    Let $f$ be a $\bw$-homogeneous polynomial on $\C^n$  of degree $d$. Given $p+q=n$, let $\ell=(p+1)d-\sum_i w_i$, and denote by $(\mathbb{M}_f)_\ell$ the linear subspace of the Milnor algebra consisting of degree $\ell$ elements.
$$h^{p,q}(K_f) = \dim_\C (\mathbb{M}_f)_\ell.
$$
    When \eqref{CYcondition} is satisfied, i.e. $K_f$ is a Calabi-Yau link, the condition reduces to $\ell = pd$.
\end{theorem}

Finally, Moriyama 
expresses the dimension of the 
moduli space of cCY structures on  a given 7-dimensional link $K_f$, in terms of the Sasakian Hodge numbers \cite{Moriyama2011}:
\begin{equation}
\label{dim: ccY}
    \dim \mathfrak{M}_{CY}(K_f, \mathcal{F}_\xi) = b^3(K_f) + h^{1,1}(K_f) - 1.
\end{equation}
In particular, the third Betti number $b^3$ is completely determined by $h^{2,1}$ and $h^{3,0}$:
\begin{equation}
    (h^{3,0}+ h^{2,1})(K_f) = \frac{1}{2}b^3(K_f),
\end{equation}
which we have computed in this work.
The remaining term is $h_S^{1,1}$, which is not calculable via Theorem \ref{th:itoh}, however may be accessible by other means.
For instance, in the study of Calabi-Yau manifolds \cite{CANDELAS1990383,Dixon:1987bg,Lerche:1989uy}, there is a well-established notion of homological mirror symmetry between Hodge numbers. 
We propose that if one could extend this to a notion of mirror symmetry among links, perhaps one could access $h^{1,1}$ for a link as being $h^{2,1}$ of the respective `mirror'. For this dataset, we have enumerated the $h^{2,1}$'s \textit{exhaustively}, so we could in principle know all the terms in  \eqref{dim: ccY} and subsequently the dimensions of the moduli space of cCY structures, at least for CY link mirror pairs.

\subsection{The Crowley-Nordström invariant on cCY $7$-manifolds}
For an arbitrary closed $7$-manifold with $\rG_{2}$-structure $(Y^{7},\varphi)$, Crowley and Nordstr\"{o}m have defined a $\mathbb{Z}/{48}\mathbb{Z}$-valued homotopy invariant  $\nu(\varphi)$, which is a combination of topological data from a compact coboundary 8-manifold with $\Spin(7)$-structure $(W^{8},\Psi)$ extending $(Y^{7},\varphi)$, in the sense that $Y = \partial W$ and $\Psi\mid_{Y} = \varphi$:
\begin{equation}
\label{eq: nu(phi) invariant}
    \nu(\varphi) := \chi(W) - 3\sigma(W) \text{ mod } 48,
\end{equation}
where $\chi$ the real Euler characteristic and $\sigma$ is the signature.

In particular ($n=3$), specialising \S \ref{sec: CY links} to real dimension $7$, a contact Calabi--Yau structure naturally induces a coclosed $\rG_2$-structure (with symmetric torsion):
\begin{proposition}[{\cite[Corollary 6.8]{Habib2015}}] \label{prop:G2estruturaCCY} \mbox{}\\
    Every cCY manifold $(K^7,\eta,\xi,J ,\Upsilon)$ is an $S^1$-bundle $\pi:K\to V$ over a Calabi--Yau 3-orbifold $(V,\omega,\Upsilon)$, with connection 1-form $\eta$ and curvature  
\begin{equation}
\label{eq:eta.curv}
    d\eta=\omega,    
\end{equation} 
    and it carries a cocalibrated $\rG_2$-structure
\begin{equation}
\label{eq:G2structure}
    \varphi 
    :=\eta\wedge \omega +\Re\Upsilon,
\end{equation} 
    with torsion $d\varphi= \omega\wedge\omega$   and Hodge dual  $ 4$-form
$$ \psi=\ast\varphi = \frac{1}{2}\omega\wedge\omega+ \eta\wedge\Im\Upsilon.$$
\end{proposition}
It therefore makes sense, for further instance, to systematically study the invariants $\nu(\varphi)$ associated to $\bw$-homogeneous polynomials $f$ defining Calabi-Yau links $K_f$.

\subsection{Weak R-equivalence of weighted polynomials}

In the light of Theorem \ref{th:itoh}, we will observe in \S\ref{sec:method}  that the Sasakian Hodge numbers and the natural CN invariant of a CY link depend only upon the Milnor algebra $\mathbb{M}_f$. The relation between the Milnor algebras of different weighted homogeneous polynomials was examined in \cite{ahmed2012}, where the notion of R-equivalence is introduced:
\begin{definition}
    Two weighted homogeneous polynomials $f,g$ on $\C^n$ are \emph{$R-$equivalent} if there exists a diffeomorphism $\psi : \C^n\mapsto\C^n$ such that $f\circ\psi = g$.
\end{definition}

Theorem 2 in \cite{ahmed2012} gives a sufficient condition for R-equivalence between two such polynomials, as quoted below:
\begin{theorem}
    Let $f,g$ be $\bw$-homogeneous polynomials on $\C^n$ of degree $d$, such that $\mathcal{J}_f = \mathcal{J}_g$; then $f$ is R-equivalent to $g$.
\end{theorem}

However, our initial empirical observations (as detailed in \S\ref{sec:explicitweakR}) suggested that any homogeneous polynomial with the same weight vector (up to permutations), and no further singularities, has the same $\ell$-degree subspaces  of the Milnor algebra, up to linear isomorphism. This then implies that their respective Sasakian Hodge numbers and CN invariants will be the same. This motivates us to propose the following definition and conjecture.

\begin{definition}
    Two weighted homogeneous polynomials $f, g$ on $\C^n$ are said to be \emph{weakly R-equivalent} if the respective $\ell$-degree linear subspaces of their Milnor algebras are isomorphic, for each $\ell$ such that $p+q=n$, as in Theorem \ref{th:itoh}.
\end{definition}

\begin{conjecture}
\label{conj: weak_R}
    Consider two weighted homogeneous polynomials $f,g$ on $\C^n$ of same degree $d$; if their weight vectors $\bw_f$ and $\bw_g$ coincide (up to permutations), then $f$ and $g$ are weakly R-equivalent.
\end{conjecture}

The Conjecture is somewhat surprising, since it encompasses cases in which the Jacobian ideals $\mathcal{J}_f$ and $\mathcal{J}_g$ are non-isomorphic, and thus $\mathbb{M}_f$ and $\mathbb{M}_g$ are not equivalent. As we will see in \S\ref{sec:method}, although certain steps in the algorithms to compute the Sasakian Hodge numbers and CN invariant involve the Gr\"obner basis, which is directly related to the Jacobian ideal, the results of these computations seem to depend only on the initial sets of weights of the $\C^\times$-action. 


\subsection{Machine Learning}

Aiming at an audience in the community of mathematics and theoretical physics, we provide a very brief introduction to neural networks, which is the architecture we use in our investigation \cite{anderson1995introduction,Ruehle:2020jrk,he2023machine}.
We begin by introducing the \emph{neuron}, the building block of any neural network. A neuron is a vertex in an oriented graph, which takes in a set of input data $\{x_{i}\}$ and produces a single numerical output $y$, by the following three steps:
\begin{enumerate}
    \item First, each input $x_{i}$ is multiplied by a weight $W_{i}$: $W_{i}x_{i}$.
    \item Next, all the weighted inputs are summed and a bias $b$ is added: $\sum_{i}W_{i}x_{i}+b$.
    \item Finally, the sum is passed through a non-linear activation function which produces an output: $\hat{y}=\text{act}(\sum_{i} W_{i}x_{i}+b)$.
\end{enumerate}
ReLU is perhaps the most standard example of a non-linear activation function, it is defined
\begin{equation}
    \text{ReLU}(x) \vcentcolon = \begin{cases} x & x>0 \\ 0 & otherwise \end{cases}
\end{equation}
and is the activation function used in this work.
A neural network is then simply a collection of neurons stratified in a series of layers, whereby the neurons in each layer are connected by edges to neurons in the previous and next layers.

The process of training a neural network starts with partitioning the dataset into \emph{training data}, from which the network will learn, and \emph{test data}, which is only used after training to evaluate the network's performance. The training process involves repeatedly calculating the `error', which is some measure of the difference between the predicted model outputs and the true known outputs for the training data. During training, these weights and biases are stochastically updated in order to reduce this error measure. Computing the error requires a choice of loss function, for regression problems typically one uses either mean absolute error (MAE) or mean squared error (MSE)
\begin{equation}
\begin{split}
\label{eq:MAE}
    \text{MAE} & = \frac{1}{N} \sum_{i=1}^{N} |y_{i} - \hat{y}_{i}|, \\
    \text{MSE} & = \frac{1}{N} \sum_{i=1}^{N} (y_{i} - \hat{y}_{i})^2,
\end{split}
\end{equation}
where $y_{i}$ and $\hat{y}_{i}$ are the true and predicted values, respectively, and $N$ is the dataset size.
The method by which we change the weights and biases to minimise the loss is called the \emph{optimisation algorithm}, the simplest of which is stochastic gradient descent (SGD). There are more advanced optimisation methods that build on SGD, of which Adam \cite{kingma2017} is a popular choice and is the particular optimisation we adopt.

For regression tasks, typical performance metrics include MAE~\eqref{eq:MAE} as well the $R^{2}$ score \eqref{eq:R2}, which is defined as the proportion of the variance in the dependent variable that is predictable from the independent variable(s): 
\begin{equation}
\label{eq:R2}
    R^{2} = 1 - \frac{\sum_{i=1}^{N}(y_{i}-\hat{y}_{i})^{2}}{\sum_{i=1}^{N}(y_{i}-\bar{y})^{2}}
    \quad\in\quad (-\infty,1],
\end{equation}
where 
\begin{equation}
    \bar{y} = \frac{1}{N} \sum_{i=1}^{N} y_{i},
\end{equation}
is the mean output. Therefore an $R^{2}$ score close to 1 means the regression model is a good fit, whereas a score close to 0 means the model is a poor fit.
Additionally, despite this being a regression problem, we introduce a classification-inspired metric: Accuracy.
We define this the be the proportion of predictions within a fixed distance from the true value over the test data, where this fixed distance is defined in terms of a bound which is $0.05$ times the range of the true values.
This also evaluates in the range $[0,1]$, where a value of 1 indicates perfect learning.

Cross-validation is a method commonly used to get an unbiased evaluation of the learning, whereby the full dataset is shuffled and then split into $k$ non-overlapping subsets. Each subset acts as the test dataset once, whilst the remaining $(k-1)$ subsets are combined to create the complement training dataset. $k$ independent identical neural network models are then each trained on one of these training datasets, then evaluated on the complementary test dataset, with the evaluation scores recorded. The mean evaluation scores, with their standard errors, are then calculated and used to measure the model performance.

Continuous datasets may be quantitatively tested for correlations via the Product Moment Correlation Coefficient (PMCC), as linear method providing first order insight about potential dataset dependencies.
This measure is defined
\begin{equation}
    \text{PMCC}(X,Y) \vcentcolon = \frac{\mathbb{E}(X-\mu_X)\mathbb{E}(Y-\mu_Y)}{\sigma_x \sigma_Y} \in [-1,1]\;,
\end{equation}
for random variables $X, Y$ with respective means $\mu_X, \mu_Y$ and standard deviations $\sigma_x, \sigma_Y$, and expectation values over the datasets taken with $\mathbb{E}(\cdot)$. 
It takes values in the range stated with perfect (anti-)correlation represented by PMCC of 1 (-1), and no correlation for PMCC of 0. 

A more interpretable supervised learning method is symbolic regression \cite{Koza1994}.
In this method a basis of functions is assembled into an expression which is fit to training data.
The basis used later in this work was $\{+,-,*,/\}$.
These methods are genetic algorithms at their core, in that a population of expressions (represented as trees) are first initialised and then evaluated on the training data with respect to a standard loss~\eqref{eq:MAE} perturbed with a parsimony term which rewards simplicity in the expression (alike regularisation techniques in traditional ML preventing overfitting).
The fittest individuals are then selected for cross-breeding between their expression trees, with subsequent mutation, to form the next generation of expressions.
This process is then iterated to convergence, providing an array of candidate expressions which one can deduce information about the true functional form from.

\section{Data Generation \& Analysis}
\label{sec:method}
As previously stated, the Calabi-Yau $3$-folds arising in the link construction are hypersurfaces in complex weighted projective space $\mathbb{P}^{4}(\mathbf{w})$. Such spaces are compact Fano manifolds (with positive curvature), constructed through identification of $\mathbb{C}^5$ with a weight vector of 5 entries. 
It was shown in \cite{CANDELAS1990383}, that the list of weight vector combinations which lead to unique weighted projective spaces whose anticanonical divisors are compact and Ricci-flat is finite, with $N=7555$ cases.

For each $\mathbb{P}^{4}(\mathbf{w})$, any hypersurface in the anticanonical divisor class can be represented as a weighted homogeneous polynomial of degree $\sum_i w_i$. 
Throughout this class, there is freedom in the choice of complex coefficients for each of the monomial terms in the hypersurface's defining polynomial equation.
Any choice of coefficients, such that the surface does not become more singular, defines a Calabi-Yau $3$-fold. 
All of these will share the same Hodge numbers, but may otherwise be topologically distinct \cite{CANDELAS1990383}.
In addition, there is redundancy between  choices of coefficient sets due to polynomial symmetries (such as coordinate transformations, coefficient normalisation, etc.), allowing multiple sets of coefficients to define the same $3$-fold.

The dataset of Calabi-Yau links considered in this work was constructed using one Calabi-Yau from each of the respective 7555 $\mathbb{P}^{4}(\mathbf{w})$'s.
In each case, the Calabi-Yau polynomial was first selected to have all monomial coefficients as $1$.
Physically, this may be interpreted as considering equivalent points on the Coulomb branches of the vacuum expectation value moduli spaces, when the Calabi-Yau manifolds are used for string compactification \cite{Feng:2000mi}. 
However 1484 out of the 7555 polynomial hypersurfaces intersected with singularities in the ambient space, leading to a higher-dimensional singularity structure on the links. 
To avoid this, for these 1484 cases other polynomials were sampled, with coefficients from $\{1,2,3,4,5\}$, until the singularity structure was exclusively the isolated singularity at the origin - as required for the link construction\footnote{Practically, the dimension of the singular locus of each Calabi-Yau polynomial was computed over a finite field of prime characteristic (101). Since this field reduction from the complex numbers cannot decrease the dimension of the singularity structure, where the dimension was 0 the polynomial was accepted. Where the observed singularity dimension was higher, a selection of other primes $(251,1993,1997)$ were used to check for bad field reduction, where the dimension was 0 in any of those cases the polynomial was accepted as the increase in observed singularity dimension was due to this bad reduction. Where the singularity dimension did not decrease to 0, the polynomial was resampled until one with singularity dimension 0 was found (each time only 1 resample was required).}.

To exemplify this process, consider the weight vector $\bw=(22,29,49,50,75)$, whose degree $d=225$ ($=\sum_i w_i$) monomial basis has 7 terms: 
$$z_1^8z_3,\; z_1^4z_2^3z_3,\; z_1z_2^7,\; z_1z_2z_3z_4z_5,\; z_2z_3^4, \;
z_4^3z_5,\; z_5^3.
$$
We thus initialise the Calabi-Yau  polynomial equation to:
\begin{align*}
    0=&\; a_1z_1^8z_3+a_2z_1^4z_2^3z_3+a_3z_1z_2^7+a_4z_1z_2z_3z_4z_5 \\&+a_5z_2z_3^4+a_6z_4^3z_5+a_7z_5^3,
\end{align*}
for the complex coefficient vector $(a_1,\dots,a_7)$. 
We set $a_1=\dots=a_7=1$, and check the singularity structure of the resulting hypersurface. 
In this case, the singular locus defined by this polynomial has dimension 0, which is the isolated singularity at the origin, there is hence no further singularity structure introduced. We therefore accept this Calabi-Yau $3$-fold, adding it to our database for topological invariant computation (no further sampling of the $a_i$ values is required). 

For 7549 of these 7555 Calabi-Yau's selected in this way,  the topological properties of the corresponding links were calculated. Namely, the Sasakian Hodge numbers $\{h^{3,0},h^{2,1}\}$, from Theorem \ref{th:itoh}, and the CN invariant, from \eqref{eq: nu(phi) invariant}.
It is worth emphasising that, since the list of weight vectors which lead to complex 3-dimensional Calabi-Yaus is finite, and since the topological invariants computed are conjectured to be identical for all Calabi-Yau polynomials with same weight vector (via Conjecture \ref{conj: weak_R}, and inspired by initial empirical observations exemplified in \S\ref{sec:explicitweakR}), the data generated for these 7555 manifolds would be \textit{exhaustive} for this link construction\footnote{In this work, the topological data was computed for 7549 of the 7555, where the 6 remaining ((1, 1, 8, 19, 28), (1, 1, 9, 21, 32), (1, 1, 11, 26, 39), (1, 1, 12, 28, 42), (1, 6, 34, 81, 122), (1, 6, 40, 93, 140)) were the most computationally expensive and whose computation is left for future work. Since they represent such a small fraction of the dataset their omission should have negligible impact on learning performance.}.

The polynomial generation and topological invariant computations were performed in \texttt{sagemath} \cite{sagemath}, with the help of \texttt{macaulay2}~\cite{M2} and \texttt{singular}~\cite{SingularLib}. 
Computation of each of the topological invariants required the respective Gr\"obner bases of the Calabi-Yau polynomials; these bases are notoriously expensive to compute with at worst doubly-exponential time complexity \cite{Dube_1990}, and taking our High-Performance Computing cluster (HPC) $\sim 100,000$ core hours.
Hence, as a side product of these computational efforts, the Gr\"obner basis for a selection of the Calabi-Yau polynomials considered (one for each possible weight vector) is provided, along with the corresponding topological quantities, on this work's respective \href{https://github.com/TomasSilva/MLcCY7.git}{\texttt{GitHub}}.

\begin{figure}[tb]
    \centering
    \includegraphics[width=0.8\linewidth]{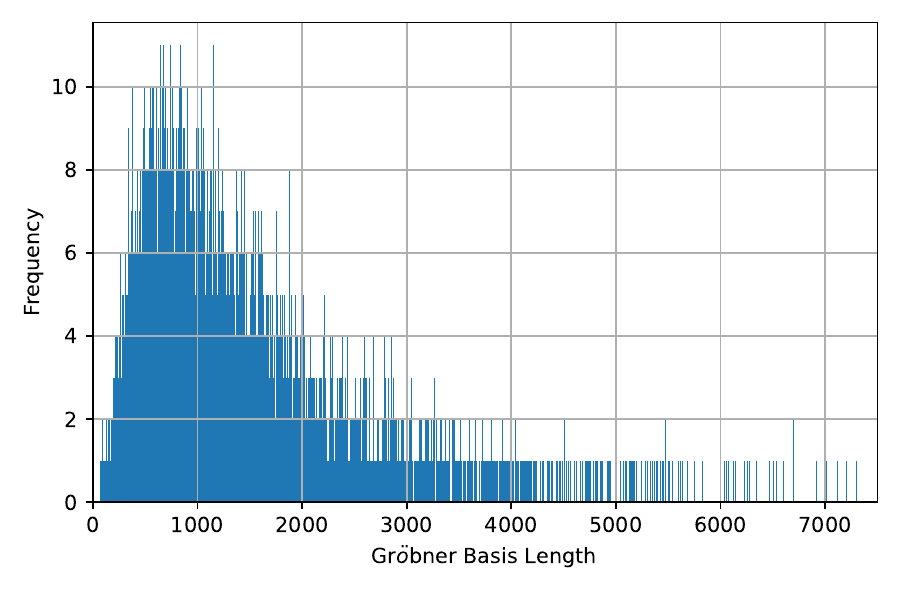} 
    \caption{Histogram of Gr\"obner basis lengths for the 7549 Calabi-Yau link constructions computed.}
    \label{grob_histogram}
\end{figure}

The distribution of the lengths of the Gr\"obner bases for 7549 out of 7555  Calabi-Yau polynomials is shown in Figure \ref{grob_histogram}. Due to the non-trivial connection between weights and basis length, and the importance of the basis length in determining whether invariant computation is even feasible, the prediction of Gr\"obner bases lengths was independently investigated, as detailed in \S\ref{grobnerml}.

\subsection{Sasakian Hodge Numbers}
\label{sec:SHN}

As outlined in \S\ref{sec:sasakian-hodge}, the computation of the Hodge numbers $h^{3,0}$ and $h^{2,1}$ associated with each weighted homogeneous polynomial is done by the algorithmic implementation of the explicit formula in Theorem~\ref{th:itoh},  from \cite{Itoh2004}:
\begin{algorithm}[ht]
    \caption{Computation of Sasakian Hodge Numbers via Theorem~\ref{th:itoh}.}
    \begin{algorithmic}[1]
        \Require $f(z_1, \dots, z_5)$, a homogeneous polynomial in $\C^5$.
        \Require $\mathbf{w}=(w_1, \dots, w_5)$, the weight vector associated with the polynomial $f$.
        \Ensure $[h^{3,0}, h^{2,1}]$, the Sasakian Hodge Numbers associated to $(f, \mathbf{w})$
        
        \State $A \gets$ $\C[z_1, \dots, z_5]$

        \State $d \gets \deg(f)$$\Comment{\text{w.r.t. }\mathbf{w}}$

        \State $J_{f} \gets \left\langle {\frac {\partial f}{\partial z_{1}}},\ldots ,{\frac {\partial f}{\partial z_{5}}}\right\rangle$$\Comment{\text{The Jacobian ideal of $f$}}$

        \State \label{algline:Groebner-computation}$K \gets \textsc{Gr\"obnerBasis}\left(\frac{A}{J_{f}}\right)$$\Comment{\text{Basis of }\mathbb{M}_f}$

        \State $h^{3,0} = \#\{x \in K : \deg(x) = 4d-\sum_i w_i\}$$\Comment{\dim (\mathbb{M}_f)_{4d-\sum w_i}}$

        \State $h^{2,1} = \#\{x \in K : \deg(x) = 3d-\sum_i w_i\}$$\Comment{\dim (\mathbb{M}_f)_{3d-\sum w_i}}$

        \State \Return $[h^{3,0}, h^{2,1}]$

    \end{algorithmic}
    \label{alg:SHN}
\end{algorithm}

Step~\ref{algline:Groebner-computation} of Algorithm~\ref{alg:SHN} corresponds to a well-known super-exponential (hard) routine, the Gr\"obner basis generation~\cite{BARDET201549, Adams_and_Philippe_Loustaunau1994-ze}. In order to perform computations within a feasible time for all the $7555$ polynomials\footnote{of which only 6 have been considered timed-out for this letter.}, we implemented a parallel version of Algorithm~\ref{alg:SHN} using \texttt{sagemath}~\cite{sagemath} and its built-in interface to \texttt{singular}~\cite{SingularLib}, and executed the job on a HPC cluster.

\begin{figure}[tb]
    \centering
    \includegraphics[width=0.8\linewidth]{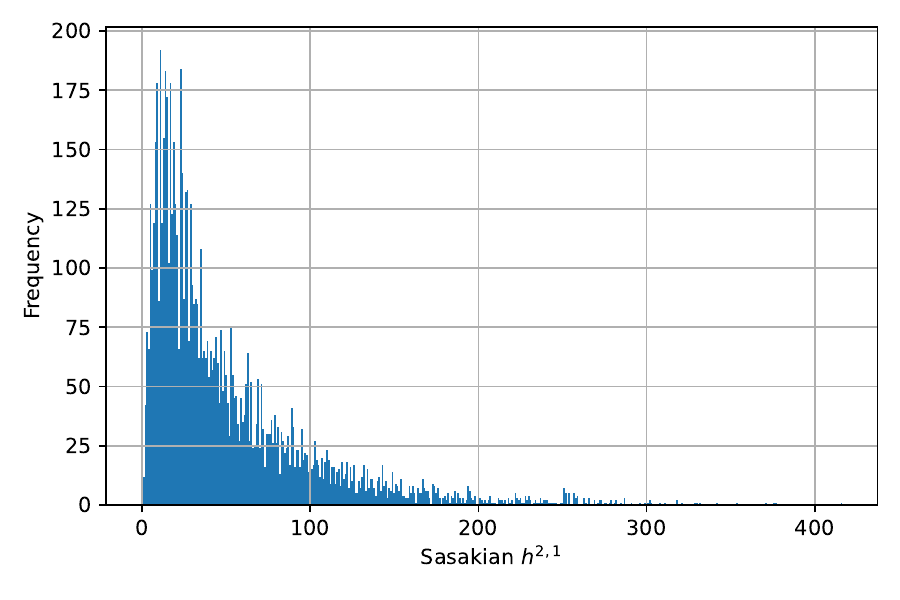} 
    \caption{Histogram of Sasakian $h^{2,1}$ values for the 7549 Calabi-Yau link constructions computed.}
    \label{h21_histogram}
\end{figure}

The Sasakian $h^{3,0}$ values for the 7549 Calabi-Yau links computed all take value 1, matching the value known for all Calabi-Yau $3$-folds, which corresponds to the unique holomorphic volume form.
The Sasakian $h^{2,1}$ values range from 1 to 416, their frequency distribution is shown in Figure \ref{h21_histogram}.
Despite the similar structure to the distribution of Gr\"obner basis lengths, we note that there is only a mild positive correlation, with PMCC $\sim 0.65$.

\begin{figure}[tb]
    \centering
    \includegraphics[width=0.8\linewidth]{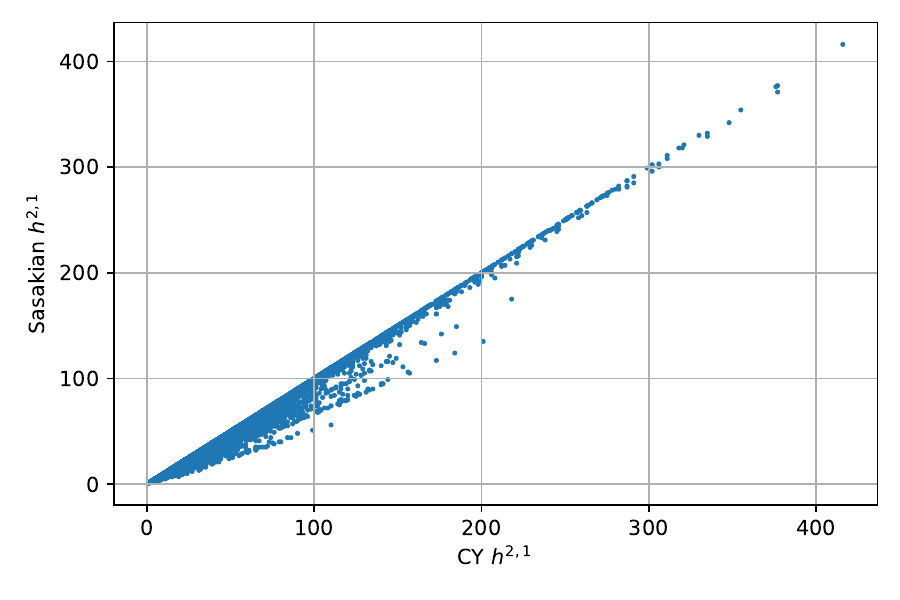} 
    \caption{Scatter graph of the Calabi-Yau complex threefold $h^{2,1}$ values against the Sasakian transverse $h^{2,1}$ values for the 7549 Calabi-Yau link constructions considered. For this data Sasakian $h^{2,1} \leq$ CY $h^{2,1}$, with 4198 cases satisfying the equality.}
    \label{h21compare_histogram}
\end{figure}

Due to the aforementioned successes of ML in predicting the Hodge numbers of these Calabi-Yaus \cite{he2017,Berman:2021mcw}, it is natural to consider  how this performance can extrapolate to Sasakian Hodge numbers, which will be  addressed in \S\ref{SHml}.
To compare directly the CY Hodge numbers with the Sasakian Hodge numbers of the links built from the same polynomials, a cross-plot of the respective $h^{2,1}$ values is given in Figure \ref{h21compare_histogram}. 
This plot shows that these topological invariants are strongly correlated (PMCC $\sim 0.99$), and the Sasakian Hodge number is bounded above by the CY Hodge number -- suggesting the following mathematical conjecture:
\begin{conjecture}
\label{conj: h21 bound}
    The Sasakian Hodge number $h_S^{2,1}$ for a Calabi-Yau link is bounded above by the Hodge number $h_{CY}^{2,1}$ of the Calabi-Yau $3$-fold built from the same $\bw$-homogeneous polynomial:
    \begin{equation}
        h_{S}^{2,1} \leq h_{CY}^{2,1}\;.
    \end{equation}
\end{conjecture}
We note that the analogous bound also technically holds for $1 = h_{S}^{3,0} \leq h_{CY}^{3,0} = 1$; as it may well be the case for other yet uncomputed Sasakian Hodge numbers.
Hence, from the successes in previous work on learning Calabi-Yau Hodge numbers \cite{he2017,Berman:2021mcw}, and this strong correlation with Sasakian Hodge numbers, the investigation into their ML prediction is well-motivated.

\subsection{Crowley-N\"ordtrom invariant}

To compute the CN invariant for polynomials in our dataset, we modify a procedure developed and described in \cite{Calvo-Andrade2020} which utilises Steenbrink's Signature theorem. Let $$
\{z^\alpha: \alpha = (\alpha_1,...,\alpha_{n+1})\in I \subset\mathbb{Z}_{>0}^{n+1}\}
$$ be a set of monomials in $\C[z_1,...,z_{n+1}]$ representing a basis over $\C$ for its Milnor algebra $\mathbb{M}_f=\frac{\C[[z_1,...,z_{n+1}]]}{\partial f/\partial z_1,...,\partial f/\partial z_{n+1}}$. For each $\alpha \in I$, define
\begin{equation}
    l(\alpha) := \sum_{i=1}^{n+1} (\alpha_i + 1)\frac{w_i}{d}.
\end{equation}
The CN invariant of a link was computed in \cite{Calvo-Andrade2020} in terms of its degree and weights, along with the signature $(\mu_-, \mu_0, \mu_+)$  of the intersection form on $H^4(V_f,\R)$:
\begin{equation}
    \nu(\varphi) = \left(\frac{d}{w_1}-1\right)...\left(\frac{d}{w_5}-1\right) - 3(\mu_+ -\mu_-) + 1
\end{equation}
Steenbrink \cite{Steenbrink1977} proved that the signature can be computed as follows:
\begin{align*}
    \mu_+ &= |\{\beta\in I: l(\beta)\notin\mathbb{Z}, \floor {l(\beta)}\in 2\mathbb{Z}\}|\\
    \mu_- &= |\{\beta\in I: l(\beta)\notin\mathbb{Z}, \floor {l(\beta)}\notin 2\mathbb{Z}\}|\\
    \mu_0 &= |\{\beta\in I: l(\beta)\in\mathbb{Z}\}|
\end{align*}
In \cite{Calvo-Andrade2020}, this procedure was originally implemented as two separate scripts, one in \texttt{singular} and one in \texttt{MATHEMATICA} \cite{Mathematica}. We improve upon this by combining those into a single \texttt{python} script. We are then able to take advantage of parallel processing, pooling, and the powerful computational resources of our HPC's to compute the CN invariant for 7549 out of the 7555 Calabi-Yau links.

\begin{figure}[tb]
    \centering
    \includegraphics[width=0.8\linewidth]{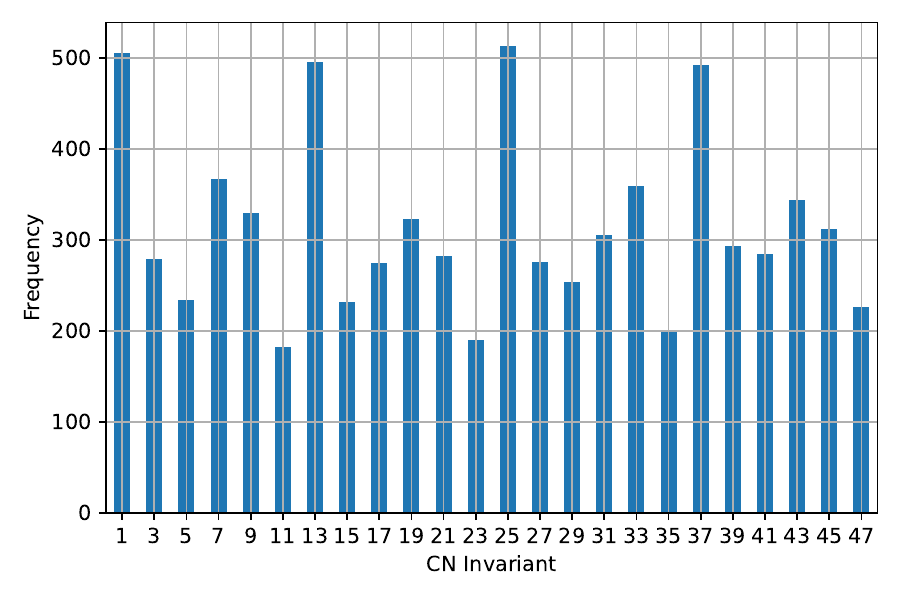} 
    \caption{Histogram of CN invariants for the 7549 Calabi-Yau link constructions computed.}
    \label{cni_histogram}
\end{figure}

The CN invariants computed fully span the range of possible values, which are odd integers from 1 to 47, cf. \cite[Proposition 3.2]{Calvo-Andrade2020}. Their frequency distribution is shown in Figure \ref{cni_histogram}, which exhibits an unexpected periodicity of 12 in the most populous invariant values ($\sim 500$). 
In particular, we note the occurrence of CN invariants 27 and 35, where previous work had not identified examples in these topological classes \cite{Calvo-Andrade2020}. 
Below we provide an explicit example of a Calabi-Yau polynomial that leads to a link in each of these classes (noting the repetition of the example considered earlier in \S\ref{sec:method}): 
\begin{align}
    & \text{CN} : 27 \,,\nonumber \\ 
    & \text{Weights}: (22, 29, 49, 50, 75) \label{eq:cCY_example}\,,\\ 
    & \text{Polynomial}: \nonumber \\
    & 0=z_1^8z_3+z_1^4z_2^3z_3+z_1z_2^7+z_1z_2z_3z_4z_5+z_2z_3^4+z_4^3z_5+z_5^3 \,, \nonumber \\
    \mbox{} \nonumber \\
    & \text{CN} : 35 \,,\nonumber \\
    & \text{Weights}: (31, 35, 36, 42, 108) \,,\\ 
    & \text{Polynomial}: \nonumber \\
    & 0=z_1^7z_2+z_1^2z_2^2z_3z_4^2+z_1z_2z_3^4z_4+z_1z_2z_3z_4z_5+z_2^6z_4+z_3^7 \nonumber \\
    & \ \ +z_3^4z_5+z_3z_5^2+z_4^6 \,. \nonumber
\end{align}

\subsection{Explicit Weak R-Equivalence}\label{sec:explicitweakR}
To corroborate the weak R-equivalence predicted in Conjecture \ref{conj: weak_R}, we show the observed behaviour for the previously selected example, as stated in \eqref{eq:cCY_example}.

In performing the checks of weak R-equivalence, we considered 10 permutations of the example weight system, and 50 polynomials per permutation (with general integer coefficients in the range (0,100), such that singular locus dimension is still 0; here, 100 was selected as it is the bound set by the prime 101 used for the coefficient ring characteristic).
In each case, the computed CN invariant was $\nu=27$, and $(h_S^{3,0},h_S^{2,1})=(1,2)$.
We observe that, while these values for the invariants were the same, the Gr\"obner basis lengths changed among different weight permutations  (but were the same for different polynomials with the same weight system permutation).
This behaviour was expected since the permutation of weights effectively amounts to a relabelling of the coordinates.

In addition to running these checks for the quoted example, the same procedure was repeated for 100 weight systems randomly selected from the database (selecting those with generally shorter polynomials, for computational efficiency), again considering 50 polynomials per weight system, and in all cases the weak R-equivalence was verified.
Code to run these checks in general scenarios is also made available in this article's \href{https://github.com/TomasSilva/MLcCY7.git}{\texttt{GitHub}}

Below are two example weight  permutations of  \eqref{eq:cCY_example}, each with two respective Calabi-Yau polynomials. These were all included in the explicit checks above and hence lead to the same topological invariants.

$$\text{Weights}: [75, 22, 49, 29, 50]\,,$$
$$0 = z_1^3+z_1z_2z_3z_4z_5+z_1z_5^3+z_2^8z_3+z_2^4z_4^3z_5+z_2z_4^7+z_3^4z_4\,,$$
$$0 = 48z_1^3+49z_1z_2z_3z_4z_5+6z_1z_5^3+71z_2^8z_3+35z_1^4z_4^3z_5+29z_2z_4^7+25z_3^4z_4\,.$$
$$\text{Weights}: [49, 22, 75, 50, 29]\,,$$
$$0 = 70z_1^4z_5+12z_1z_2^8+39z_1z_2z_3z_4z_5+90z_2^4z_4z_5^3+95z_2z_5^7+49z_3^3+11z_3z_4^3\,,$$
$$0 = 30z_1^4z_5+22z_1z_2^8+23z_1z_2z_3z_4z_5+59z_2^4z_4z_5^3+90z_2z_5^7+38z_3^3+7z_3z_4^3\,.$$

\section{Machine Learning}
\label{sec:results}

In order to investigate the efficacy of ML techniques to learn topological invariants of this dataset of Calabi-Yau links, NNs were chosen as the prototypical tool from supervised learning.
Since the output invariants take a large range of values in each case, the NNs were set up for a regression-style problem.
The NNs used had the same architecture in each case. They had neuron layer sizes of $(16,32,16)$, ReLU activation, and were trained on a MSE loss using an Adam optimiser.
These layer sizes and the other hyperparameters were set after some heuristic tuning for the Gr\"obner basis ML, then used for the other investigations also for consistency.
Each NN hence amounts to a map of the form:
\begin{equation}\label{eq:nnlayers}
    \mathbb{R}^5 \xrightarrow[]{f_1} \mathbb{R}^{16} \xrightarrow[]{f_2} \mathbb{R}^{32} \xrightarrow[]{f_3} \mathbb{R}^{16} \xrightarrow[]{f_4} \mathbb{R}^1\;,
\end{equation}
such that each $f_i$ acts via linear then non-linear action as $f(\textbf{x}) = \text{ReLU}(\textbf{W}\cdot\textbf{x} + \textbf{b})$.

In each case, the regression NNs were trained on 5 different partitions of the dataset into 80:20 train:test splits in accordance with cross-validation, to provide statistical error on the metrics used to assess learning performance.
The NNs were implemented in \texttt{python} with the use of \texttt{scikit-learn} \cite{scikit-learn}.

These NNs were trained to predict, from an input of the weight vectors, the respective Calabi-Yau link Gr\"obner basis length, Sasakian Hodge number $h^{2,1}$, and the CN invariant.
The first two of these investigations are detailed in the subsequent sections, whilst the final investigation is briefly detailed here since this architecture could not learn to predict the CN invariant sufficiently well.
Performance had $R^2$ value of $\sim 0.004$. 
Even reducing the problem to a binary classification between the two most populous classes ($\nu = 1$ and $\nu = 25$) did not lead to accuracies much above 0.5, indicating no significant learning and highlighting the highly non-trivial dependency of this invariant on the input polynomial and weight data, despite our computations showing the invariant was only dependent on the weight data in accordance with Conjecture~\eqref{conj: weak_R}.

\subsection{Grobner Basis Length}\label{grobnerml}
By a significant margin, the computational bottleneck (in terms of both time and memory) of the invariant calculations was the generation of the Gr\"obner basis.
Many initial runs failed from memory overload in this step for specific Calabi-Yau links.
Through trial and error, we diverted our resources to allocate more computational power to the harder cases with larger Gr\"obner bases. However, it would have been substantially more efficient if we had approximately accurate predictions of which links would have required more computational power.

This problem again suits itself to ML, where a simple regression model can provide quick estimates of Gr\"obner basis length, and thus guide the subsequent allocation of computational resources. 
Previously, ML methods have been used to help optimise specific steps of the Gr\"obner basis algorithm \cite{peifer2020learning, Mojsilovi__2023}, and decide when a Gr\"obner basis would be useful \cite{Huang_2016}.
Yet, only recently has ML been used to predict Gr\"obner basis properties directly, eg. in \cite{jamshidi2023predicting} predictions of Gr\"obner basis length reached $R^2$ scores $\sim 0.4$ for binomial ideals of 2 terms.

For our specific construction the polynomial ideals have much more than 2 terms (with the longest polynomial having 680), and 5 terms per ideal for each partial derivative of the polynomial.
Therefore it is prudent to investigate the performance of regression NNs in predicting Gr\"obner basis length for our Calabi-Yau links. 
These NNs used the same hyperparameters as previously specified in \eqref{eq:nnlayers}, with neuron layer sizes (16,32,16), ReLU activation, training on a MSE loss with an Adam optimiser on a 80:20 train:test data split. 
Performance measures for the learning were:
\begin{equation}
\begin{split}
    R^2 & =  0.964 \pm 0.002 \,,\\
    \text{MAE} & =  \ \ \ 122 \pm 2\,,\\
    \text{Accuracy} & = 0.947 \pm 0.005\,.
\end{split}
\end{equation}
showing excellent performance, especially when comparing the MAE score to the Gr\"obner basis range of values, as shown in Figure \ref{grob_histogram}.
We also note that whilst the number of monomial terms in a polynomial does correlate with Gr\"obner basis length, the correlation is not strong (PMCC $\sim 0.66$), and hence ML is much more useful in estimating the length of a Gr\"obner basis.
These results corroborate the motivation to use ML in improving computational efficiency, particularly here in the application of calculating topological invariants; where perhaps the architectures are using that the data represents similar geometries to aid learning in this case. 
Inspired by this success, the authors hope to extend, in future work, the application of ML to learning more properties of Gr\"obner bases, as well as the basis elements directly.

\subsection{Sasakian Hodge Numbers}\label{SHml}
Extending work where NNs have showed surprising success in predicting Hodge numbers of Calabi-Yau manifolds, we now investigate their success in predicting the Sasakian Hodge number $h^{2,1}$, the computation of which was described in \S\ref{sec:SHN}.

Using the same NN regressor architecture previously described in \eqref{eq:nnlayers}, with neuron layer sizes (16,32,16), ReLU activation, training on a MSE loss with an Adam optimiser on a 80:20 train:test data split. 
The $h^{2,1}$ values for the 7549 Calabi-Yau links computed were learnt with performance measures:
\begin{equation}\label{eq:NNHresults}
\begin{split}
    R^2 & =  0.969 \pm 0.003 \,,\\
    \text{MAE} & =  \ \ 5.53 \pm 0.22\,,\\
    \text{Accuracy} & = 0.967 \pm 0.004\,.
\end{split}
\end{equation}
These results are equivalently strong and exemplify the efficacy of ML methods in predicting more subtle topological parameters. 

Comparing the predicted $h^{2,1}$ outputs of a trained NN to the $h_{CY}^{2,1}$ values of the input weight system (instead of the intended $h_S^{2,1}$ values) produces lower performance scores across the dataset: $R^2 = 0.915$, MAE $= 9.88$, Accuracy $= 0.912$.
This reassures that the NN's are learning the Sasakian topology instead of the correlated (and well learnt in previous work) Calabi-Yau base topology.

\subsection{Sasakian Hodge Symbolic Regression}\label{SHsr}

Motivated by the highly accurate regression results, one is led to expect these NN functions to be approximating a true relation between the Sasakian Hodge numbers and the weights used to define the Calabi-Yau links.
In spirit, this could be a similar phenomenon to what was observed for the Poincar\'e polynomial, cf. \cite{Batyrev:2020ych}.

To distil some mathematical insight about the function space of appropriate approximations for Sasakian Hodge numbers from weights, as equivalently and independently probed by the NN architectures, we implement techniques of symbolic regression using \texttt{PySR}~\cite{cranmer2023interpretable} to provide interpretable relations between inputs and outputs, which should help guide any investigation of this direct relation bypassing the Milnor algebra Gr\"obner basis computation. The scripts used for the symbolic regression analysis are also available at this work's \href{https://github.com/TomasSilva/MLcCY7.git}{\texttt{GitHub}}. The (highest performing) equation on the independent test data (10\% of the dataset) proposed by \texttt{PySR} to model Sasakian $h^{2,1}$ is presented in \eqref{eq:SymReg_eq}, achieving $R^2\approx0.99$ and $\MAE\approx 2.6$, exceeding the scores of the NN in \eqref{eq:NNHresults}.
\begin{align}
\begin{split}
    \label{eq:SymReg_eq}
    h^{2,1}_\texttt{PySR}(w_0, \dots, w_4)= &\frac{14.91 w_{1} \left(w_{0} w_{4} + w_{3} \left(w_{0} + w_{3}\right)\right)}{w_{0} w_{1} w_{2} w_{3}}\\
    &+\frac{10.02 w_{2} w_{3} \left(w_{0} + w_{4} + 0.77\right)}{w_{0} w_{1} w_{2} w_{3}}\;,
\end{split}
\end{align}
where this expression may shed some light on the structure of the true function.

\begin{figure*}[!t]
    \centering
    \begin{minipage}{\columnwidth}
        \includegraphics[scale=0.45]{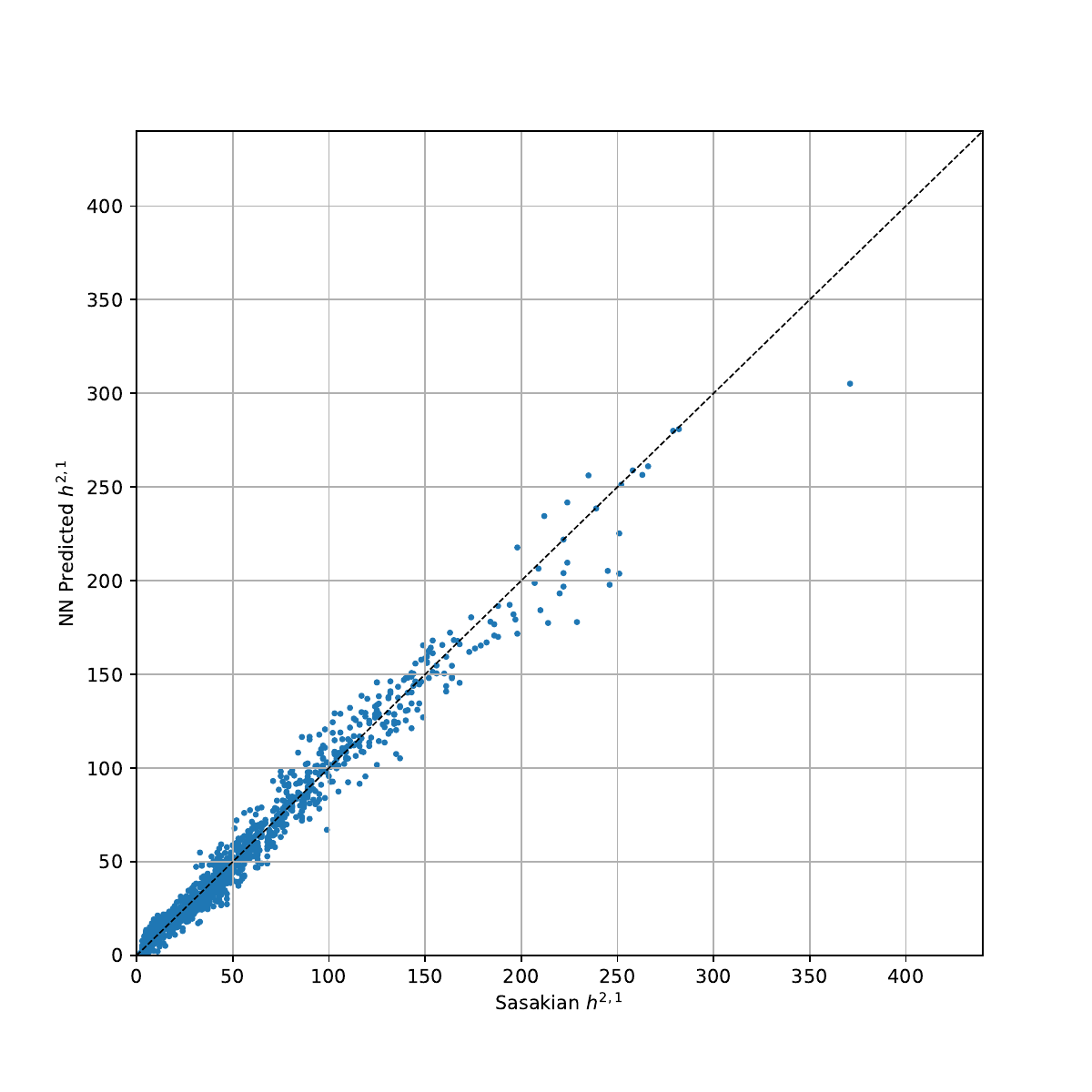}
        \label{fig:trueprd_NN}
        \vspace{-20pt}
        \par $\qquad\qquad\qquad\qquad\qquad$ \footnotesize{(a) Neural Network}
    \end{minipage}
    \hfill
    \begin{minipage}{\columnwidth}
        \includegraphics[scale=0.45]{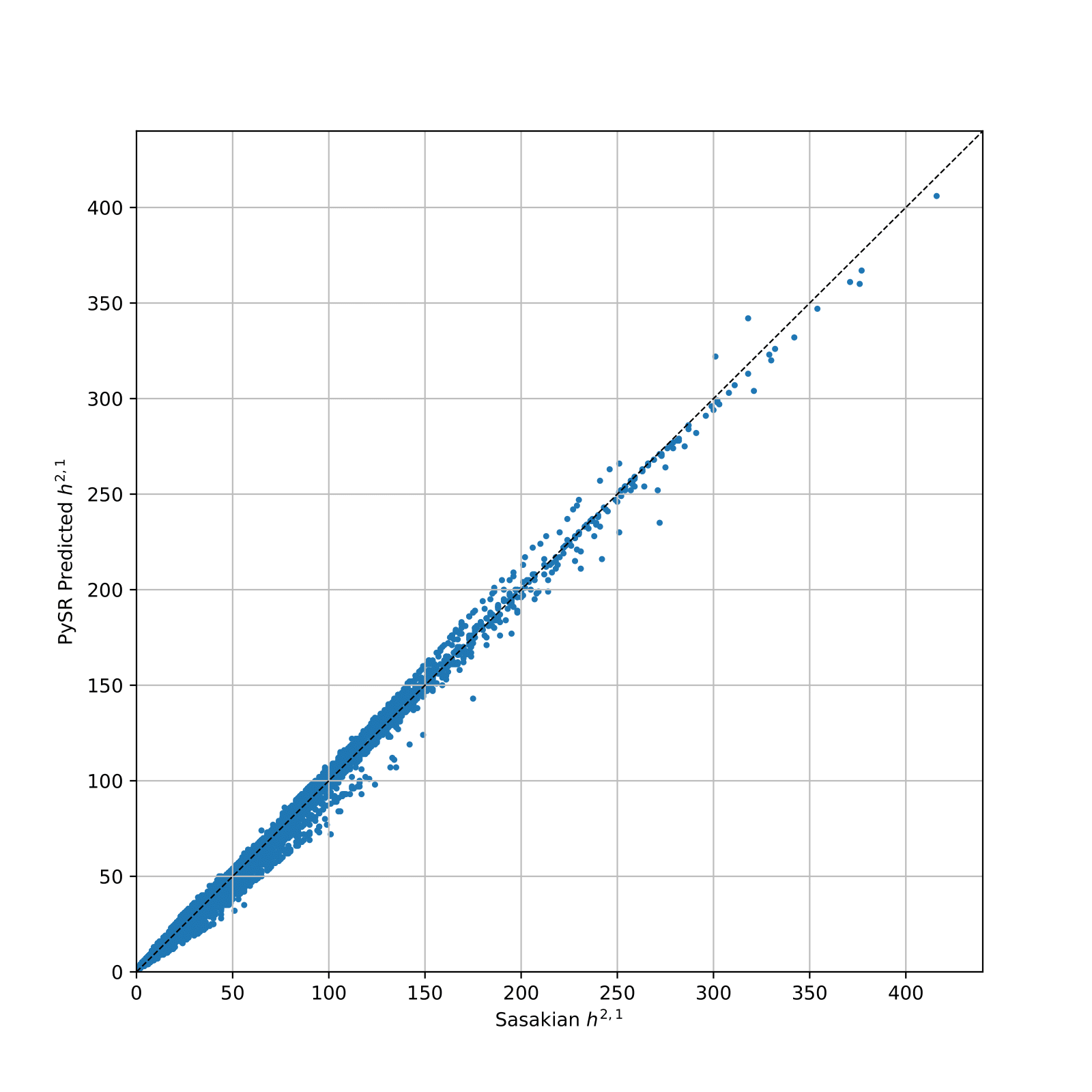}
        \label{fig:trueprd_SR}
        \vspace{-20pt}
        \par $\qquad\qquad\qquad\qquad\quad\ $ \footnotesize{(b) Symbolic Regression}
    \end{minipage}
    \vspace{-5pt}
    \caption{ML architecture predictions of the $h_S^{2,1}$ values, against the true values, for the 7549 Calabi-Yau link constructions considered, from: (a) a trained NN; and (b) the symbolic regression best model of equation \eqref{eq:SymReg_eq}.}
    \label{fig:truepred}
\end{figure*}

The $h_S^{2,1}$ predicted value from \eqref{eq:SymReg_eq} for each of the Calabi-Yau links considered is plotted against the true values in Figure \ref{fig:truepred}.
Additionally the plot shows the equivalent predictions for a trained NN, where predictions are less accurate, particularly at higher invariant values.

\subsection{Predicting the Remaining Invariants}
\label{inv_predict}

As mentioned in \S\ref{sec:method}, in this work the topological invariants have been computed for all but 6 of the 7555 weight systems, where HPC walltime has been successively exceeded due to the enormous computational cost.
Practically the computations failed at walltimes of $1000$ core hours using $>200$GB of RAM.
However, as demonstrated in the previous sections, machine learning architectures have shown great success in modelling the prediction of these topological invariants. It therefore makes sense to apply instances of these very trained models to predict the invariant values for the remaining (currently computing) six weight systems.
The prediction results are shown in Table \ref{tab:final6_predictions}.

\begin{table}[!t]
\centering
\begin{tabular}{|c|c|cc|}
\hline
Invariant           & GBL  & \multicolumn{2}{c|}{$h_S^{2,1}$} \\ \hline
Architecture        & NN   & \multicolumn{1}{c|}{NN}    & SR  \\ \hline
(1, 1, 8, 19, 28)   & 1531  & \multicolumn{1}{c|}{257}   & 338   \\ \hline
(1, 1, 9, 21, 32)   & 1623  & \multicolumn{1}{c|}{272}   & 377   \\ \hline
(1, 1, 11, 26, 39)  & 1807  & \multicolumn{1}{c|}{293}   & 447   \\ \hline
(1, 1, 12, 28, 42)  & 1879  & \multicolumn{1}{c|}{299}   & 476   \\ \hline
(1, 6, 34, 81, 122) & 3488 & \multicolumn{1}{c|}{243}   & 243   \\ \hline
(1, 6, 40, 93, 140) & 3770 & \multicolumn{1}{c|}{265}   & 272   \\ \hline
\end{tabular}
\caption{ML architecture predictions of the well-learnt CY link properties for the remaining 6 weight systems. The Gr\"obner basis length (GBL) and the Sasakian $h_S^{2,1}$ values are predicted by the trained NN architectures or the symbolic regression (SR) equation.}
\label{tab:final6_predictions}
\end{table}

Due to the excessive computation time for these remaining 6 weight systems, we expect their Gr\"obner basis lengths to make up the tail end of the respective distribution in Figure \ref{grob_histogram}. 
The extrapolation capabilities of the NN models were briefly tested by training the same NN architecture to predict Gr\"obner basis length from input weight system for 95\% of the dataset pairs with shortest lengths, then testing on the remaining 5\% ($378$ bases).
The NN had lower performance scores of: $R^2 = 0.594$, MAE $= 452$, Accuracy $= 0.484$; highlighting the importance of interpretation of statistical methods for out-of-distribution predictions.
However the predictions do well identify these test weight systems as having long bases, with the minimum predicted basis length on the test data being 2391, far exceeding the mean of the full dataset (1290). 
Whilst the NN models do predict large lengths for the remaining weight systems' Gr\"obner bases in Table \ref{tab:final6_predictions}, with values significantly above the average, they do not predict values exceeding the current highest value of 7299.
These predictions are high, as is useful for our implementation but less likely to be exact.

Conversely, the Sasakian Hodge number values have been shown to only loosely correlate with the Gr\"obner basis length, hence there is only loose intuition that these Hodge numbers should be high.
Our Conjecture \ref{conj: h21 bound} produces some bound on the Sasakian $h^{2,1}$ values, set by the respective Calabi-Yau $h^{2,1}$ values: 348, 387, 462, 491, 246, 275.
Interestingly, all predictions from both the NN and the symbolic regression satisfy this bound as well.
Furthermore, the better-performing symbolic regression model of \eqref{eq:SymReg_eq} predicts Hodge numbers very close to this bound, which is sensible behaviour, as demonstrated by Figure \ref{h21compare_histogram}.

Despite the predictions between the NN and symbolic regressor in some cases being quite different, pragmatically mathematicians and physicists are often most interested in discerning when Hodge numbers are particularly low, such as 0 or 1, as indications of phenomena such as exactness, rigidity or unobstruction. 
Hence, these learnt models can very quickly provide confidence in whether this is the case, by making predictions which can test conjectures or corroborate theoretical expectations.

Our own experienced difficulty in computing the topological invariants for these final 6 remaining cases illustrates  the potential of these machine learning models.
Where direct computation is not feasible, ML methods can provide predictions for quantities of interest (such as our CY link topological invariants) with statistical confidence, providing invaluable insight to guide refinement of the computation, and further the progress of academic research.

\section{Conclusion}
\label{sec:conclusion}

In this work, real 7-dimensional Calabi-Yau links were constructed for a complex 3-dimensional Calabi-Yau manifold coming from 7549 of the 7555 complex 4-dimensional weighted projective spaces that admit them.
It was observed, and conjectured, that any Calabi-Yau hypersurface with the correct singularity properties lead to the same Sasakian Hodge number $h^{2,1}$ and CN invariant; which (when computed for the final 6 evasive cases) will produce an exhaustive list of these invariants from the Calabi-Yau link construction. 

The datasets of these invariants were statistically analysed, and NN regressors were used to successfully predict the respective Sasakian $h^{2,1}$ values from the ambient $\mathbb{P}^4(\textbf{w})$ weights alone.
The same architectures were not successful in predicting the CN invariant, but did show surprising success in predicting the length of the Gr\"obner basis from the weights.
These regressors can hence be used to streamline future computation of Gr\"obner basis by informing on efficient computational resource allocation.

The exhaustive list of Calabi-Yau link data, as well as the \texttt{python} scripts used for their analysis and ML are made available at this work's corresponding \href{https://github.com/TomasSilva/MLcCY7.git}{\texttt{GitHub}}.

Avenues for future work include studying whether we can obtain Sasakian structure for more general links made through general toric varieties rather than weighted projective spaces. 
As well as application of ML techniques to study other Gr\"obner basis properties to further streamline their computation.

\appendix
\section{Sasakian Hodge decomposition}
\label{sec: Appendix}

  We establish some notation and elementary facts about the complexified tangent bundle, which are largely adapted from \cite{Biswas2010}.  The contact structure splits the tangent bundle as $TM=  B\oplus N_\xi$, where $B=\ker(\eta)$ and $N_\xi$ is the real line bundle spanned by the Reeb field $\xi$. The transverse complex structure $\Phi$  satisfies $(\Phi\vert_{ B })^{2}=-1$,  so the  eigenvalues of the complexified operator $\Phi^{\C}$ are $\pm\sqrt{-1}$. The complexification $B_\C:=B\otimes_\R\C$ splits   as 
  $B_\C = B^{1,0}\oplus  B ^{0,1}$, so we obtain a decomposition of direct sum of vector bundles 
$$
  \Lambda^{d}(B_\C)^\ast  
  =\bigoplus_{i=0}^{d} (B^{i,d-i})^\ast.
$$
  This induces the decomposition of vector bundles
\[
  \Omega^{d}(M)=\left(\bigoplus_{i=0}^{d} \Omega_{B }^{i,d-i}(M) \right) \oplus\left(\eta\wedge\left( \bigoplus_{j=0}^{d-1}\Omega_{B }^{j,d-j-1}(M)\right)\right), 
\]
  where $\Omega^{p,q}_{B}(M):=\Gamma(M,\Lambda^p(B_\C)^\ast\otimes\Lambda^q(B_\C)^\ast)$.

Let us briefly describe the transverse complex geometry on a Sasakian manifold $(M,\eta,\xi,g,\Phi)$, following \cite{Portilla2023} and references therein. Relatively to the Reeb foliation, the usual Hodge star induces a \emph{transverse  Hodge star} operator $ \ast_T\colon \Omega_B^k(M) \to \Omega_B^{m-(k+1)}(M)  $ by the formula 
\begin{equation}
\label{eq:transverse star}
    \ast_T(\beta)  =(-1)^{m-1-k}\ast(\beta\wedge\eta).
\end{equation}
Both operators are compatible, in the sense that 
\begin{equation}
\label{eq:star transverse star}
    \ast\alpha= \ast_T\alpha\wedge \eta,
    \quad\forall
    \alpha\in\Omega^\bullet_B(M).
\end{equation} 

 Now, acting on $(p,q)$-forms, we have well-defined operators 
\begin{equation}
    \label{eq: partial_B}
\begin{array}{lll}
  \partial_B :\Omega^{p,q}_B\to\Omega^{p+1,q}_B   &\text{and } &
  \bar{\partial}_B :\Omega^{p,q}_B\to\Omega^{p,q+1}_B.
 \end{array}
\end{equation}
In terms of the transverse Hodge star \eqref{eq:transverse star}, the operators  \eqref{eq: partial_B} have adjoints
\begin{equation}
\label{eq:adjoin partial_B}
    \begin{array}{ll }
    \partial_B^\ast\colon       &\Omega^{p,q}_B\to\Omega^{p-1,q}_B  \\  [4pt]
    &  \partial_B^\ast:= -\ast_T \partial_B\ast_T   
 \end{array}
 \qandq
 \begin{array}{ll }
    \bar{\partial}_B^\ast\colon &\Omega^{p,q}_B\to\Omega^{p,q-1}_B\\ [4pt] 
    &\bar{\partial}_B^\ast:=-\ast_T\bar{\partial}_B\ast_T.
 \end{array}
\end{equation}

  An inner product   on $\Omega^\bullet(M)$ is defined  by: 
\begin{equation}
\label{eq:innerProduct}
    (\alpha,\beta)_M=\displaystyle\int_M\alpha\wedge\ast\beta=\int_M \alpha\wedge\ast_T\beta \wedge\eta.
\end{equation} 
  
  Defining  $\Delta_{\bar{\partial}_B} :=\bar{\partial}_B^\ast\bar{\partial}_B +\bar{\partial}_B\bar{\partial}_B^\ast$, $H^{p,q}:=\ker \Delta_{\bar{\partial}_B}\subset \Omega_B^{p,q} $, the space $H^k$ of harmonic $k$-forms on $M$ has the following decomposition~\cite{Tanaka1975, Itoh2004}:
\begin{equation}
\label{eq:HodgeDec}
    H^k = \bigoplus_{p+q=k} H^{p,q}, \qwithq
    \bar{H}^{p,q} = H^{q,p}, \quad k \leq n.
\end{equation}
\newpage
\section*{Acknowledgements}

\noindent\underline{YHH} is supported by a Leverhulme Trust grant and an STFC grant ST/J00037X/2. \underline{E.~Heyes} is supported by City, University of London and the States of Jersey. \underline{E.~Hirst} is supported by Pierre Andurand. \underline{HSE} is supported by FAPESP grants 2018/21391-1, 2020/09838-0 and 2021/ 04065-6, the Brazilian National Council for Scientific and Technological Development (CNPq) grant \mbox{307217/2017-5}, and UK Royal Society grant NMG$\backslash$R1$\backslash$191068. \underline{TSRS} is supported by São Paulo Research Foundation (FAPESP) grant 2022/09891-4.

\par This research utilised Queen Mary's Apocrita HPC facility, supported by QMUL Research-IT \cite{apocrita}; and computational resources of the ``\textit{Centro Nacional de Processamento de Alto Desempenho em São Paulo} (CENAPAD-SP)".


\appendix

\bibliographystyle{elsarticle-num} 
\bibliography{references.bib}

\begin{thebibliography}{10}
\expandafter\ifx\csname url\endcsname\relax
  \def\url#1{\texttt{#1}}\fi
\expandafter\ifx\csname urlprefix\endcsname\relax\def\urlprefix{URL }\fi
\expandafter\ifx\csname href\endcsname\relax
  \def\href#1#2{#2} \def\path#1{#1}\fi

\bibitem{he2017}
Y.-H. He, {Deep-Learning the Landscape} (6 2017).
\newblock \href {http://arxiv.org/abs/1706.02714} {\path{arXiv:1706.02714}}.

\bibitem{ruehle2017}
F.~Ruehle, Evolving neural networks with genetic algorithms to study the string
  landscape, Journal of High Energy Physics 2017~(8) (aug 2017).
\newblock \href {https://doi.org/10.1007/jhep08(2017)038}
  {\path{doi:10.1007/jhep08(2017)038}}.

\bibitem{carifio2017}
J.~Carifio, J.~Halverson, D.~Krioukov, B.~D. Nelson, Machine learning in the
  string landscape, Journal of High Energy Physics 2017~(9) (sep 2017).
\newblock \href {https://doi.org/10.1007/jhep09(2017)157}
  {\path{doi:10.1007/jhep09(2017)157}}.

\bibitem{krefl2017}
D.~Krefl, R.-K. Seong, Machine learning of {Calabi-Yau} volumes, Physical
  Review D 96~(6) (sep 2017).
\newblock \href {https://doi.org/10.1103/physrevd.96.066014}
  {\path{doi:10.1103/physrevd.96.066014}}.

\bibitem{Bull:2019cij}
K.~Bull, Y.-H. He, V.~Jejjala, C.~Mishra, {Getting {CICY} High}, Phys. Lett. B
  795 (2019) 700--706.
\newblock \href {http://arxiv.org/abs/1903.03113} {\path{arXiv:1903.03113}},
  \href {https://doi.org/10.1016/j.physletb.2019.06.067}
  {\path{doi:10.1016/j.physletb.2019.06.067}}.

\bibitem{He:2020lbz}
Y.-H. He, A.~Lukas, {Machine Learning {Calabi-Yau} Four-folds}, Phys. Lett. B
  815 (2021) 136139.
\newblock \href {http://arxiv.org/abs/2009.02544} {\path{arXiv:2009.02544}},
  \href {https://doi.org/10.1016/j.physletb.2021.136139}
  {\path{doi:10.1016/j.physletb.2021.136139}}.

\bibitem{Erbin_2021}
H.~Erbin, R.~Finotello, Machine learning for complete intersection {Calabi-Yau}
  manifolds: A methodological study, Physical Review D 103~(12) (jun 2021).
\newblock \href {https://doi.org/10.1103/physrevd.103.126014}
  {\path{doi:10.1103/physrevd.103.126014}}.

\bibitem{Berman:2021mcw}
D.~S. Berman, Y.-H. He, E.~Hirst, {Machine learning {Calabi-Yau}
  hypersurfaces}, Phys. Rev. D 105~(6) (2022) 066002.
\newblock \href {http://arxiv.org/abs/2112.06350} {\path{arXiv:2112.06350}},
  \href {https://doi.org/10.1103/PhysRevD.105.066002}
  {\path{doi:10.1103/PhysRevD.105.066002}}.

\bibitem{pmlr-v197-aslan23a}
B.~Aslan, D.~Platt, D.~Sheard,
  \href{https://proceedings.mlr.press/v197/aslan23a.html}{Group invariant
  machine learning by fundamental domain projections}, in: S.~Sanborn,
  C.~Shewmake, S.~Azeglio, A.~Di~Bernardo, N.~Miolane (Eds.), Proceedings of
  the 1st NeurIPS Workshop on Symmetry and Geometry in Neural Representations,
  Vol. 197 of Proceedings of Machine Learning Research, PMLR, 2023, pp.
  181--218.
\newline\urlprefix\url{https://proceedings.mlr.press/v197/aslan23a.html}

\bibitem{ashmore2020}
A.~Ashmore, Y.-H. He, B.~A. Ovrut, Machine learning {Calab-Yau} metrics,
  Fortschritte der Physik 68~(9) (2020) 2000068.
\newblock \href {https://doi.org/10.1002/prop.202000068}
  {\path{doi:10.1002/prop.202000068}}.

\bibitem{anderson2020}
L.~B. Anderson, M.~Gerdes, J.~Gray, S.~Krippendorf, N.~Raghuram, F.~Ruehle,
  {{Moduli-dependent Calabi-Yau and SU(3)-structure metrics from Machine
  Learning}}, JHEP 05 (2021) 013.
\newblock \href {http://arxiv.org/abs/2012.04656} {\path{arXiv:2012.04656}},
  \href {https://doi.org/10.1007/JHEP05(2021)013}
  {\path{doi:10.1007/JHEP05(2021)013}}.

\bibitem{jejjala2020}
V.~Jejjala, D.~K. Mayorga~Pena, C.~Mishra, {Neural network approximations for
  Calabi-Yau metrics}, JHEP 08 (2022) 105.
\newblock \href {http://arxiv.org/abs/2012.15821} {\path{arXiv:2012.15821}},
  \href {https://doi.org/10.1007/JHEP08(2022)105}
  {\path{doi:10.1007/JHEP08(2022)105}}.

\bibitem{douglas2020}
M.~R. Douglas, S.~Lakshminarasimhan, Y.~Qi, {Numerical {Calabi-Yau} metrics
  from holomorphic networks} (12 2020).
\newblock \href {http://arxiv.org/abs/2012.04797} {\path{arXiv:2012.04797}}.

\bibitem{larfors2021}
M.~Larfors, A.~Lukas, F.~Ruehle, R.~Schneider, {Learning Size and Shape of
  Calabi-Yau Spaces} (11 2021).
\newblock \href {http://arxiv.org/abs/2111.01436} {\path{arXiv:2111.01436}}.

\bibitem{klaewer2019}
D.~Klaewer, L.~Schlechter, Machine learning line bundle cohomologies of
  hypersurfaces in toric varieties, Physics Letters B 789 (2019) 438--443.
\newblock \href {https://doi.org/10.1016/j.physletb.2019.01.002}
  {\path{doi:10.1016/j.physletb.2019.01.002}}.

\bibitem{berglund2023new}
P.~Berglund, Y.-H. He, E.~Heyes, E.~Hirst, V.~Jejjala, A.~Lukas, {{New
  Calabi-Yau Manifolds from Genetic Algorithms}} (6 2023).
\newblock \href {http://arxiv.org/abs/2306.06159} {\path{arXiv:2306.06159}}.

\bibitem{Manko:2022zfz}
M.~Manko, {An Upper Bound on the Critical Volume in a Class of Toric
  {Sasaki-Einstein} Manifolds} (9 2022).
\newblock \href {http://arxiv.org/abs/2209.14029} {\path{arXiv:2209.14029}}.

\bibitem{He:2020eva}
Y.-H. He, E.~Hirst, T.~Peterken, {Machine-learning dessins d'enfants:
  explorations via modular and {Seiberg-Witten} curves}, J. Phys. A 54~(7)
  (2021) 075401.
\newblock \href {http://arxiv.org/abs/2004.05218} {\path{arXiv:2004.05218}},
  \href {https://doi.org/10.1088/1751-8121/abbc4f}
  {\path{doi:10.1088/1751-8121/abbc4f}}.

\bibitem{Bao:2020nbi}
J.~Bao, S.~Franco, Y.-H. He, E.~Hirst, G.~Musiker, Y.~Xiao, {Quiver Mutations,
  {S}eiberg Duality and Machine Learning}, Phys. Rev. D 102~(8) (2020) 086013.
\newblock \href {http://arxiv.org/abs/2006.10783} {\path{arXiv:2006.10783}},
  \href {https://doi.org/10.1103/PhysRevD.102.086013}
  {\path{doi:10.1103/PhysRevD.102.086013}}.

\bibitem{Bao:2021vxt}
J.~Bao, O.~Foda, Y.-H. He, E.~Hirst, J.~Read, Y.~Xiao, F.~Yagi, {Dessins
  d'enfants, {Seiberg-Witten} curves and conformal blocks}, JHEP 05 (2021) 065.
\newblock \href {http://arxiv.org/abs/2101.08843} {\path{arXiv:2101.08843}},
  \href {https://doi.org/10.1007/JHEP05(2021)065}
  {\path{doi:10.1007/JHEP05(2021)065}}.

\bibitem{Bao:2021auj}
J.~Bao, Y.-H. He, E.~Hirst, J.~Hofscheier, A.~Kasprzyk, S.~Majumder, {{H}ilbert
  series, machine learning, and applications to physics}, Phys. Lett. B 827
  (2022) 136966.
\newblock \href {http://arxiv.org/abs/2103.13436} {\path{arXiv:2103.13436}},
  \href {https://doi.org/10.1016/j.physletb.2022.136966}
  {\path{doi:10.1016/j.physletb.2022.136966}}.

\bibitem{Bao:2021olg}
J.~Bao, Y.-H. He, E.~Hirst, {Neurons on Amoebae}, J. Symb. Comput. 116 (2022)
  1--38.
\newblock \href {http://arxiv.org/abs/2106.03695} {\path{arXiv:2106.03695}},
  \href {https://doi.org/10.1016/j.jsc.2022.08.021}
  {\path{doi:10.1016/j.jsc.2022.08.021}}.

\bibitem{Bao:2021ofk}
J.~Bao, Y.-H. He, E.~Hirst, J.~Hofscheier, A.~Kasprzyk, S.~Majumder, {Polytopes
  and Machine Learning} (9 2021).
\newblock \href {http://arxiv.org/abs/2109.09602} {\path{arXiv:2109.09602}}.

\bibitem{AriasTamargo:2022qgb}
G.~Arias-Tamargo, Y.-H. He, E.~Heyes, E.~Hirst, D.~Rodriguez-Gomez, {Brain webs
  for brane webs}, Phys. Lett. B 833 (2022) 137376.
\newblock \href {http://arxiv.org/abs/2202.05845} {\path{arXiv:2202.05845}},
  \href {https://doi.org/10.1016/j.physletb.2022.137376}
  {\path{doi:10.1016/j.physletb.2022.137376}}.

\bibitem{Dechant:2022ccf}
P.-P. Dechant, Y.-H. He, E.~Heyes, E.~Hirst, {Cluster Algebras: Network Science
  and Machine Learning} (3 2022).
\newblock \href {http://arxiv.org/abs/2203.13847} {\path{arXiv:2203.13847}}.

\bibitem{Chen:2022jwd}
S.~Chen, Y.-H. He, E.~Hirst, A.~Nestor, A.~Zahabi, {{M}ahler Measuring the
  Genetic Code of Amoebae} (12 2022).
\newblock \href {http://arxiv.org/abs/2212.06553} {\path{arXiv:2212.06553}}.

\bibitem{Cheung:2022itk}
M.-W. Cheung, P.-P. Dechant, Y.-H. He, E.~Heyes, E.~Hirst, J.-R. Li,
  {Clustering Cluster Algebras with Clusters} (12 2022).
\newblock \href {http://arxiv.org/abs/2212.09771} {\path{arXiv:2212.09771}}.

\bibitem{Ashmore:2023ajy}
A.~Ashmore, Y.-H. He, E.~Heyes, B.~A. Ovrut, {{N}umerical spectra of the
  {L}aplacian for line bundles on {Calabi-Yau} hypersurfaces} (5 2023).
\newblock \href {http://arxiv.org/abs/2305.08901} {\path{arXiv:2305.08901}}.

\bibitem{bao2022}
J.~Bao, Y.-H. He, E.~Heyes, E.~Hirst, {Machine Learning Algebraic Geometry for
  Physics} (4 2022).
\newblock \href {http://arxiv.org/abs/2204.10334} {\path{arXiv:2204.10334}}.

\bibitem{he2023machine}
Y.-H. He, E.~Heyes, E.~Hirst, {Machine Learning in Physics and Geometry} (3
  2023).
\newblock \href {http://arxiv.org/abs/2303.12626} {\path{arXiv:2303.12626}}.

\bibitem{Acharya2001}
B.~Acharya, E.~Witten, Chiral fermions from manifolds of $\rm{G}_2$ holonomy
  (2001).
\newblock \href {http://arxiv.org/abs/arXiv:hep-th/0109152v1}
  {\path{arXiv:arXiv:hep-th/0109152v1}}.

\bibitem{Acharya2004}
B.~Acharya, S.~Gukov, {M} theory and singularities of exceptional holonomy
  manifolds, Physics Reports 392~(3) (2004) 121--189.

\bibitem{CANDELAS1990383}
P.~Candelas, M.~Lynker, R.~Schimmrigk,
  \href{https://www.sciencedirect.com/science/article/pii/055032139090185G}{{Calabi-Yau}
  manifolds in weighted {P4}}, Nuclear Physics B 341~(2) (1990) 383--402.
\newblock \href {https://doi.org/https://doi.org/10.1016/0550-3213(90)90185-G}
  {\path{doi:https://doi.org/10.1016/0550-3213(90)90185-G}}.
\newline\urlprefix\url{https://www.sciencedirect.com/science/article/pii/055032139090185G}

\bibitem{delaOssa:2014lma}
X.~de~la Ossa, M.~Larfors, E.~E. Svanes, {{Exploring $SU(3)$ structure moduli
  spaces with integrable $G_2$ structures}}, Adv. Theor. Math. Phys. 19 (2015)
  837--903.
\newblock \href {http://arxiv.org/abs/1409.7539} {\path{arXiv:1409.7539}},
  \href {https://doi.org/10.4310/ATMP.2015.v19.n4.a5}
  {\path{doi:10.4310/ATMP.2015.v19.n4.a5}}.

\bibitem{delaOssa:2017pqy}
X.~de~la Ossa, M.~Larfors, E.~E. Svanes, {The Infinitesimal Moduli Space of
  Heterotic {$G_2$} Systems}, Commun. Math. Phys. 360~(2) (2018) 727--775.
\newblock \href {http://arxiv.org/abs/1704.08717} {\path{arXiv:1704.08717}},
  \href {https://doi.org/10.1007/s00220-017-3013-8}
  {\path{doi:10.1007/s00220-017-3013-8}}.

\bibitem{delaOssa:2017gjq}
X.~de~la Ossa, M.~Larfors, E.~E. Svanes, {Restrictions of Heterotic {$G_2$}
  Structures and Instanton Connections}, in: {Nigel Hitchin's 70th Birthday
  Conference}, 2017.
\newblock \href {http://arxiv.org/abs/1709.06974} {\path{arXiv:1709.06974}}.

\bibitem{Lotay2023}
J.~Lotay, H.~S. Earp, The heterotic $\rm{G}_2$ system on contact {Calabi--Yau}
  $7$-manifolds, Transactions of the American Mathematical Society, Series B
  10~(26) (2023) 907--943.
\newblock \href {https://doi.org/10.1090/btran/129}
  {\path{doi:10.1090/btran/129}}.

\bibitem{tomassini2007contact}
A.~Tomassini, L.~Vezzoni, Contact {Calabi-Yau} manifolds and special
  {L}egendrian submanifolds (2007).
\newblock \href {http://arxiv.org/abs/math/0612232}
  {\path{arXiv:math/0612232}}.

\bibitem{Habib2015}
G.~Habib, L.~Vezzoni, Some remarks on {C}alabi-{Y}au and hyper-{K}\"ahler
  foliations, Differential Geom. Appl. 41 (2015) 12--32.
\newblock \href {https://doi.org/10.1016/j.difgeo.2015.03.006}
  {\path{doi:10.1016/j.difgeo.2015.03.006}}.

\bibitem{Calvo-Andrade2020}
O.~Calvo-Andrade, L.~Rodr{\'i}guez, H.~N. S{\'a}~Earp,
  \href{https://www.ems-ph.org/journals/show_pdf.php?issn=0213-2230&vol=36&iss=6&rank=6}{{Gauge
  theory and $G_2$-geometry on Calabi-Yau links}}, Rev. Mat. Iberoam. 36~(6)
  (2020) 1753--1778.
\newblock \href {http://arxiv.org/abs/1606.09271 [math.DG]}
  {\path{arXiv:1606.09271 [math.DG]}}.
\newline\urlprefix\url{https://www.ems-ph.org/journals/show_pdf.php?issn=0213-2230&vol=36&iss=6&rank=6}

\bibitem{Itoh2004}
M.~Itoh,
  \href{https://www.jstage.jst.go.jp/article/kyushujm/58/1/58_1_121/_article/-char/ja/}{{S}asakian
  manifolds, {H}odge decomposition and {M}ilnor algebras}, Kyushu J. Math 58
  (2004) 121--140.
\newline\urlprefix\url{https://www.jstage.jst.go.jp/article/kyushujm/58/1/58_1_121/_article/-char/ja/}

\bibitem{CNInvariant}
D.~Crowley, J.~Nordstr\"om, New invariants of {$G2$}-structures, Geometry \&
  Topology 19 (2015) 2949--2992.

\bibitem{Vafa:1989xc}
C.~Vafa, {{S}tring Vacua and Orbifoldized {L-G} Models}, Mod. Phys. Lett. A 4
  (1989) 1169.
\newblock \href {https://doi.org/10.1142/S0217732389001350}
  {\path{doi:10.1142/S0217732389001350}}.

\bibitem{github}
D.~Aggarwal, Y.-H. He, E.~Heyes, E.~Hirst, H.~N.~S. Earp, T.~S.~R. Silva,
  \href{https://github.com/TomasSilva/MLcCY7}{{MLcCY7}} (2023).
\newline\urlprefix\url{https://github.com/TomasSilva/MLcCY7}

\bibitem{Portilla2023}
L.~E. Portilla, H.~N. S{\'a}~Earp, {Instantons on Sasakian 7-manifolds}, The
  Quarterly Journal of Mathematics (2023).
\newblock \href
  {http://arxiv.org/abs/https://academic.oup.com/qjmath/advance-article-pdf/doi/10.1093/qmath/haad011/49696708/haad011.pdf}
  {\path{arXiv:https://academic.oup.com/qjmath/advance-article-pdf/doi/10.1093/qmath/haad011/49696708/haad011.pdf}},
  \href {https://doi.org/10.1093/qmath/haad011}
  {\path{doi:10.1093/qmath/haad011}}.

\bibitem{Boyer2008}
C.~P. Boyer, K.~Galicki, Sasakian geometry, Oxford Mathematical Monographs,
  Oxford University Press, Oxford, 2008.

\bibitem{Milnor1969}
J.~Milnor, {Singular points of complex hypersurfaces} (1969).

\bibitem{Steenbrink1977}
J.~Steenbrink, Intersection form for quasi-homogeneous singularities,
  Compositio Math. 34~(2) (1977) 211--223.

\bibitem{Steenbrink1983}
J.~H.~M. Steenbrink, {Mixed Hodge structures associated with isolated
  singularities} (1983).
\newblock \href {https://doi.org/10.1090/PSPUM/040.2/713277}
  {\path{doi:10.1090/PSPUM/040.2/713277}}.

\bibitem{Moriyama2011}
T.~Moriyama, {The moduli space of transverse Calabi-Yau structures on foliated
  manifolds} (2009).
\newblock \href {https://doi.org/10.18910/12555} {\path{doi:10.18910/12555}}.

\bibitem{Dixon:1987bg}
L.~J. Dixon, {Some world sheet properties of Superstring compactifications, on
  Orbifolds and otherwise}, in: {Summer Workshop in High-energy Physics and
  Cosmology}, 1987.

\bibitem{Lerche:1989uy}
W.~Lerche, C.~Vafa, N.~P. Warner, {{C}hiral Rings in {N=2} Superconformal
  Theories}, Nucl. Phys. B 324 (1989) 427--474.
\newblock \href {https://doi.org/10.1016/0550-3213(89)90474-4}
  {\path{doi:10.1016/0550-3213(89)90474-4}}.

\bibitem{ahmed2012}
I.~Ahmed, {W}eighted homogeneous polynomials with isomorphic {M}ilnor algebras,
  Journal of Prime Research in Mathematics (2012).

\bibitem{anderson1995introduction}
J.~A. Anderson, An introduction to neural networks, MIT press, 1995.

\bibitem{Ruehle:2020jrk}
F.~Ruehle, {Data science applications to string theory}, Phys. Rept. 839 (2020)
  1--117.
\newblock \href {https://doi.org/10.1016/j.physrep.2019.09.005}
  {\path{doi:10.1016/j.physrep.2019.09.005}}.

\bibitem{kingma2017}
D.~P. Kingma, J.~Ba, {Adam}: A method for stochastic optimization (2017).
\newblock \href {http://arxiv.org/abs/1412.6980} {\path{arXiv:1412.6980}}.

\bibitem{Koza1994}
J.~R. Koza, Genetic programming as a means for programming computers by natural
  selection, Statistics and Computing 4~(2) (1994) 87--112.
\newblock \href {https://doi.org/10.1007/BF00175355}
  {\path{doi:10.1007/BF00175355}}.

\bibitem{Feng:2000mi}
B.~Feng, A.~Hanany, Y.-H. He, {D-brane gauge theories from toric singularities
  and toric duality}, Nucl. Phys. B 595 (2001) 165--200.
\newblock \href {http://arxiv.org/abs/hep-th/0003085}
  {\path{arXiv:hep-th/0003085}}, \href
  {https://doi.org/10.1016/S0550-3213(00)00699-4}
  {\path{doi:10.1016/S0550-3213(00)00699-4}}.

\bibitem{sagemath}
T.~S. Developers, W.~Stein, D.~Joyner, D.~Kohel, J.~Cremona, B.~Eröcal,
  \href{http://www.sagemath.org}{{SageMath, version 9.0}} (2020).
\newline\urlprefix\url{http://www.sagemath.org}

\bibitem{M2}
D.~R. Grayson, M.~E. Stillman, {Macaulay2, a software system for research in
  algebraic geometry}, Available at \url{http://www.math.uiuc.edu/Macaulay2/}
  (1992).

\bibitem{SingularLib}
W.~Decker, G.-M. Greuel, G.~Pfister, H.~Sch\"onemann, {{\sc Singular} {4-3-0}
  --- {A} computer algebra system for polynomial computations},
  \url{http://www.singular.uni-kl.de} (2022).

\bibitem{Dube_1990}
T.~W. Dub\'{e}, The structure of polynomial ideals and {G}röbner bases, SIAM
  Journal on Computing 19~(4) (1990) 750--773.
\newblock \href {https://doi.org/10.1137/0219053} {\path{doi:10.1137/0219053}}.

\bibitem{BARDET201549}
M.~Bardet, J.-C. Faugère, B.~Salvy,
  \href{https://www.sciencedirect.com/science/article/pii/S0747717114000935}{{On
  the complexity of the \texttt{F5} Gr\"obner basis algorithm}}, Journal of
  Symbolic Computation 70 (2015) 49--70.
\newblock \href {https://doi.org/https://doi.org/10.1016/j.jsc.2014.09.025}
  {\path{doi:https://doi.org/10.1016/j.jsc.2014.09.025}}.
\newline\urlprefix\url{https://www.sciencedirect.com/science/article/pii/S0747717114000935}

\bibitem{Adams_and_Philippe_Loustaunau1994-ze}
A.~William, P.~Loustaunau, {An introduction to Grobner bases}, American
  Mathematical, 1994.

\bibitem{Mathematica}
W.~R. Inc., \href{https://www.wolfram.com/mathematica}{{M}athematica, {V}ersion
  13.1}, champaign, IL, 2022.
\newline\urlprefix\url{https://www.wolfram.com/mathematica}

\bibitem{scikit-learn}
F.~Pedregosa, G.~Varoquaux, A.~Gramfort, V.~Michel, B.~Thirion, O.~Grisel,
  M.~Blondel, P.~Prettenhofer, R.~Weiss, V.~Dubourg, J.~Vanderplas, A.~Passos,
  D.~Cournapeau, M.~Brucher, M.~Perrot, E.~Duchesnay, {Scikit-learn: Machine
  Learning in Python}, Journal of Machine Learning Research 12 (2011)
  2825--2830.

\bibitem{peifer2020learning}
D.~Peifer, M.~Stillman, D.~Halpern-Leistner, Learning selection strategies in
  {B}uchberger's algorithm (2020).
\newblock \href {http://arxiv.org/abs/2005.01917} {\path{arXiv:2005.01917}}.

\bibitem{Mojsilovi__2023}
J.~Mojsilovi{\'{c}}, D.~Peifer, S.~Petrovi{\'{c}}, Learning a performance
  metric of {B}uchberger's algorithm, Involve, a Journal of Mathematics 16~(2)
  (2023) 227--248.
\newblock \href {https://doi.org/10.2140/involve.2023.16.227}
  {\path{doi:10.2140/involve.2023.16.227}}.

\bibitem{Huang_2016}
Z.~Huang, M.~England, J.~H. Davenport, L.~C. Paulson, Using machine learning to
  decide when to precondition cylindrical algebraic decomposition with
  {G}roebner bases (sep 2016).
\newblock \href {https://doi.org/10.1109/synasc.2016.020}
  {\path{doi:10.1109/synasc.2016.020}}.

\bibitem{jamshidi2023predicting}
S.~Jamshidi, E.~Kang, S.~Petrović, Predicting the cardinality of a reduced
  {G}r\"obner basis (2023).
\newblock \href {http://arxiv.org/abs/2302.05364} {\path{arXiv:2302.05364}}.

\bibitem{Batyrev:2020ych}
V.~V. Batyrev, {On the stringy Hodge numbers of mirrors of quasi-smooth
  {Calabi-Yau} hypersurfaces} (6 2020).
\newblock \href {http://arxiv.org/abs/2006.15825} {\path{arXiv:2006.15825}}.

\bibitem{cranmer2023interpretable}
M.~Cranmer, {Interpretable Machine Learning for Science with PySR and
  SymbolicRegression.jl} (2023).
\newblock \href {http://arxiv.org/abs/2305.01582} {\path{arXiv:2305.01582}}.

\bibitem{Biswas2010}
I.~Biswas, G.~Schumacher, Vector bundles on {S}asakian manifolds, Adv. Theor.
  Math. Phys. 14~(2) (2010) 541--562.

\bibitem{Tanaka1975}
N.~Tanaka, \href{http://hdl.handle.net/2433/84914}{A differential geometric
  study on strongly pseudo-convex manifolds}, Vol.~9, Kinokuniya, 1975.
\newline\urlprefix\url{http://hdl.handle.net/2433/84914}

\bibitem{apocrita}
T.~King, S.~Butcher, L.~Zalewski, {Apocrita - High Performance Computing
  Cluster for Queen Mary University of London} (Mar. 2017).
\newblock \href {https://doi.org/10.5281/zenodo.438045}
  {\path{doi:10.5281/zenodo.438045}}.

\end{thebibliography}

\end{document}